\documentclass[12pt] {article}
\usepackage{amsmath,amsthm}
\usepackage{amssymb,latexsym}
\usepackage{graphicx}
\usepackage{amssymb}
\usepackage{amsfonts}
\usepackage{epic,eepic}
\usepackage{amscd}
\usepackage{amscd}
\title{New structures on valuations and applications.}
\date{}
\author{Semyon Alesker \footnote{Partially supported by ISF grant 701/08.}
\\  { \normalsize Department of Mathematics, Tel Aviv University, Ramat Aviv}
 \\  { \normalsize 69978 Tel Aviv,
Israel }
\\ {\normalsize e-mail: semyon@post.tau.ac.il}}
\def\One{{1\hskip-2.5pt{\rm l}}}
\def\eps{\varepsilon}
\def\alp{\alpha}
\def\ome{\omega}
\def\Ome{\Omega}
\def\lam{\lambda}
\def\Lam{\Lambda}

\def\to{\rightarrow}

\def\grc{{}\!^ {\textbf{C}} Gr}

\def\qed { Q.E.D. }
\def\pt{\partial}

\def\RR{\mathbb{R}}
\def\CC{\mathbb{C}}

\def\PP{\mathbb{P}}

\def\FF{\mathbb{F}}

\newtheorem{theorem}{Theorem}[subsection]
\newtheorem{corollary}[theorem]{Corollary}

\newtheorem{proposition}[theorem]{Proposition}

\theoremstyle{definition}
\newtheorem{example}[theorem]{Example}
\newtheorem{definition}[theorem]{Definition}
\newtheorem{remark}[theorem]{Remark}
\theoremstyle{conjecture}

  \def\cf{{\cal F}}
  
 \def\ck{{\cal K}} 
\def\cm{{\cal M}}  
\def\cp{{\cal P}} 
\def\car{{\cal R}}
  
\def\cv{{\cal V}} \def\cw{{\cal W}} 
 
\def\inj{\hookrightarrow}

\numberwithin{equation}{section}
\begin{document}
\maketitle \setcounter{section}{-1}

\tableofcontents

\section{Introduction}\label{S:introduction}
The theory of valuations on convex sets is a classical part of
convexity with traditionally strong relations to integral geometry.
Approximately during the last 15 years there was a considerable
progress in the valuations theory and its applications to integral
geometry. The progress is both conceptual and technical: several new
structures on valuations have been discovered, new classification
results of various special classes of valuations have been obtained,
the tools used in the valuations theory and the relations of it to
other parts of mathematics have become much more diverse (besides
convexity and integral geometry, one can mention representation
theory, geometric measure theory, elements of contact geometry and
complex and quaternionic analysis). This progress in the valuations
theory has led to new developments in integral geometry,
particularly in Hermitian spaces. Some of the new structures turned
out to encode in an elegant and useful way an important integral
geometric information: for example the product on valuations encodes
somehow the principal kinematic formulas in various spaces.

%************************************
%The progress is both conceptual and technical, and it led to new
%advances in integral geometry, especially in Hermitian spaces.

Quite recently, generalizations of the classical theory of
valuations on convex sets to the context of manifolds were
initiated; this development extends the applicability of the
valuations theory beyond affine spaces, and also covers a broader
scope of integral geometric problems. In particular, the theory of
valuations on manifolds provides a common point of view on three
classical and previously unrelated directions of integral geometry:
Crofton style integral geometry dealing with integral geometric and
differential geometric invariants of sets and their intersections
with and projections to lower dimensional subspaces; Gelfand style
integral geometry dealing with the Radon transform on smooth
functions on various spaces; and less classical but still well known
the Radon transform with respect to the Euler characteristic on
constructible functions.

The relations between the valuations theory and the Crofton style
integral geometry have been known since the works of Blaschke and
especially Hadwiger, but the new developments have enriched the both
subjects, and in fact more progress is expected. The relations of
the valuations theory to the two other types of integral geometry
are new.

Besides new notions, theorems, and applications, these recent
developments contain a fair amount of new intuition on the subject.
However when one tries to make this intuition formally precise, the
clarity of basic ideas is often lost among numerous technical
details; moreover in a few cases this formalization has not been
done yet. Here in several places I take the opportunity to use the
not very formal format of lecture notes to explain the new intuition
in a heuristic way, leaving the technicalities aside. Nevertheless I
clearly separate formal rigorous statements from such heuristic
discussions.

%While the classical valuation theory was applied to integral
%geometry of Euclidean spaces, new developments led to a progress in
%other geometries, most notably in Hermitian integral geometry.

The goal of my and Joe Fu's lectures is to provide an introduction
to these modern developments. These two sets of lectures complement
each other. My lectures concentrate mostly on the valuations theory
itself and provide a general background for Fu's lectures. In my
lectures the discussion of the relations between the valuations
theory and integral geometry is usually relatively brief, and its
goal is to give simple illustrations of general notions. The
important exceptions are Sections \ref{Ss:radon} and
\ref{Ss:inversion} where new integral geometric results are
discussed, namely a Radon type transform on valuations. Much more
thorough discussion of applications to Crofton style integral
geometry, especially in Hermitian spaces, will be given in Fu's
lectures.

\hfill

My lectures consist of two main parts. The first part discusses the
theory of valuations on convex sets, and the second part discusses
its recent generalizations to manifolds. The theory of valuations on
convex sets is a very classical and much studied area. In these
lectures I mention only several facts from these classical
developments which are necessary for our purposes; I refer to the
surveys \cite{mcmullen-schneider}, \cite{mcmullen-survey} for
further details and history.

These lectures contain almost no proofs. I tried to give the
necessary definitions and list the main properties and sometimes
present constructions of the principal objects and some intuition
behind. Among important new operations on valuations are product,
convolution, Fourier type transform, pull-back, push-forward, and
the Radon type transform on valuations; all of them are relevant to
integral geometry and are discussed in these notes.

\hfill

Several interesting recent developments in the valuations theory are
not discussed here. The main omissions are a series of
investigations by M. Ludwig with collaborators of valuations with
weaker assumptions on continuity and various symmetries (see e.g.
\cite{ludwig-reitzner-99}, \cite{ludwig-duke-03},
\cite{ludwig-reitzner-annals}) and convex bodies valued valuations
(see e.g. \cite{ludwig-adv-02}, \cite{ludwig-ajm-06},
\cite{schneider-schuster}). Particularly let me mention the
surprising Ludwig-Reitzner characterization
\cite{ludwig-reitzner-99} of the affine surface area as the only (up
to the Euler characteristic, volume, and a non-negative
multiplicative factor) example of upper semi-continuous convex
valuation invariant under all affine volume preserving
transformations.

\hfill

{\bf Acknowledgements.} These are notes of my lectures to be given
at Centre de Recerca Matem\`atica during the Advanced Course on
Integral Geometry and Valuation Theory; I thank this institution and
the organizers of the course E. Gallego, X. Gual, G. Solanes, and E.
Teufel, for the invitation to give these lectures. I thank A. Bernig
for his remarks on the first version of the notes, and F. Schuster
for a very careful reading of them and numerous remarks.

\section{Translation invariant valuations on convex
sets.}\label{S:translation-invariant}

\subsection{Definitions.}\label{Ss:definitions-translation}
Let $V$ be a finite dimensional vector space of dimension $n$.
Throughout these notes we will denote by $\ck(V)$ the family of all
convex compact non-empty subsets of $V$.
\begin{definition}\label{D:valuation}
A complex valued functional
$$\phi\colon \ck(V)\to \CC$$
is called a valuation if
$$\phi(A\cup B)=\phi(A)+\phi(B)-\phi(A\cap B)$$
whenever $A,B,A\cup B\in \ck(V)$.
\end{definition}
\begin{remark}\label{R:val-rem}
In Section \ref{S:valuations-on-manifolds} we will introduce a
different but closely related notion of valuation on a smooth
manifold. To avoid abuse of terminology, we will sometimes call
valuations on convex sets in from Definition \ref{D:valuation} by
{\itshape convex valuations}, though this is not a traditional
terminology. But when it does not lead to abuse of terminology, we
will call them just valuations. In fact all valuations from Section
\ref{S:translation-invariant} will be convex, while from Section
\ref{S:valuations-on-manifolds} will not unless otherwise stated.
\end{remark}

\begin{example}\label{Ex:valuations}
(1) Any $\CC$-valued measure on $V$ is a convex valuation. In
particular the Lebesgue measure is.

(2) The Euler characteristic $\chi$ defined by $\chi(K)=1$ for any
$K\in \ck(V)$, is a convex valuation.

(3) Let $\phi$ be a convex valuation. Let $C\in \ck(V)$ be fixed.
Define
$$\psi(K):=\phi(K+C).$$
Then $\psi$ is a convex valuation. (Here $K+C:=\{k+c|k\in K,\, c\in
C\}$ is the Minkowski sum.) Indeed $(A\cup B)+C=(A+C)\cup (B+C)$,
and if $A,B,A\cup B\in \ck(V)$ then
$$(A\cap B)+C=(A+ C)\cap(B+ C).$$
\end{example}

Let us define a very important class of continuous convex
valuations. Fix a Euclidean metric on $V$. The Hausdorff distance on
$\ck(V)$ is defined by
$$dist_H(A,B):=\inf\{\eps >0|\, A\subset (B)_\eps, B\subset
(A)_\eps\}$$ where $(A)_\eps$ denotes the $\eps$-neighborhood of $A$
in the Euclidean metric. It is well known (see e.g.
\cite{schneider-book}) that $\ck(V)$ equipped with $dist_H$ is a
metric locally compact space in which any closed bounded set is
compact.
\begin{definition}
A convex valuation $\phi\colon \ck(V)\to \CC$ is called continuous
if $\phi$ is continuous in the Hausdorff metric.
\end{definition}
Readily this notion of continuity of a valuation is independent of
the choice of a Euclidean metric on $V$.

\begin{definition}
A convex valuation $\phi\colon \ck(V)\to \CC$ is called translation
invariant if
$$\phi(K+x)=\phi(K)\mbox{ for any } K\in \ck(V),\, x\in V.$$
\end{definition}

The linear space of translation invariant continuous convex
valuations will be denoted by $Val(V)$. Equipped with the topology
of uniform convergence on compact subsets of $\ck(V)$, $Val(V)$ is a
Fr\'echet space. In fact it follows from McMullen's decomposition
(Corollary \ref{C:mcmullen-decomposition} below) that $Val(V)$ with
this topology is a Banach space with a norm is given by
$$||\phi||:=\sup_{K\subset D}|\phi(K)|,$$
where $D\subset V$ is the unit Euclidean ball for some auxiliary
Euclidean metric.

\subsection{McMullen's theorem and mixed volumes.}\label{Ss:mcmullen}
The following result due to McMullen \cite{mcmullen-77} is very
important.
\begin{theorem}\label{T:mcmullen}
Let $\phi\colon \ck(V)\to \CC$ be a translation invariant continuous
convex valuation. Then for any convex compact sets $A_1,\dots,A_s\in
\ck(V)$ the function $$f(\lam_1,\dots,\lam_s):=\phi(\lam_1
A_1+\dots+ \lam_s A_s)$$ with $\lam_1,\dots,\lam_s\geq 0$ is a
polynomial of degree at most $n=\dim V$.
\end{theorem}
The special case $s=1$ is already non-trivial and important. It
means that for $\lam\geq 0$
$$\phi(\lam K)=\phi_0(K)+\lam \phi_1(K)+\dots+\lam^n\phi_n(K).$$
It is easy to see that the coefficients $\phi_0,\phi_1,\dots,\phi_n$
are also continuous translation invariant convex valuations.
Moreover $\phi_i$ is homogeneous of degree $i$ (or $i$-homogeneous
for brevity). By definition, a valuation $\psi$ is called
$i$-homogeneous if for any $K\in\ck(V),\, \lam\geq 0$ one has
$$\phi(\lam K)=\lam^i\phi(K).$$
Let us denote by $Val_i(V)$ the subspace in $Val(V)$ of
$i$-homogeneous valuations. We immediately get the following
corollary:
\begin{corollary}[McMullen's decomposition]\label{C:mcmullen-decomposition}
$$Val(V)=\oplus_{i=0}^n Val_i(V).$$
\end{corollary}

\begin{remark}\label{R:classification-extremal-cases}
Clearly $Val_0(V)$ is one dimensional and is spanned by the Euler
characteristic. Actually $Val_n(V)$ is also one dimensional and is
spanned by a Lebesgue measure; this fact is not obvious and was
proved by Hadwiger \cite{hadwiger-book}.
\end{remark}

\hfill

Let us now recall the definition of (Minkowski's) mixed volumes
which provide interesting examples of translation invariant
continuous convex valuations. Fix a Lebesgue measure $vol$ on $V$.
For any $n$ tuple of convex compact sets $A_1,\dots, A_n$ consider
the function
$$f(\lam_1,\dots,\lam_n)=vol(\lam_1A_1+\dots+\lam_nA_n).$$
This is a homogeneous polynomial in $\lam_i\geq 0$ of degree $n$. Of
course, this fact follows from McMullen's theorem \ref{T:mcmullen}
and $n$-homogeneity of the volume, though originally it was
discovered much earlier by Minkowski, and in this particular case
there is a simpler proof (see e.g. \cite{schneider-book}).
\begin{definition}\label{D:mixed-vol}
The coefficient of the monomial $\lam_1\dots,\lam_n$ in the
polynomial $f(\lam_1,\dots,\lam_n) $ divided by $n!$ is called the
mixed volume of $A_1,\dots,A_n$ and is denoted by
$V(A_1,\dots,A_n)$.
\end{definition}
The normalization of the coefficient is chosen in such a way that
$V(A,\dots,A)=vol(A)$. Mixed volumes have a number of interesting
properties, in particular they satisfy the Aleksandrov-Fenchel
inequality \cite{schneider-book}. The property relevant for us
however is the valuation property. Fix $1\leq s\leq n-1$ and an
$s$-tuple of convex compact sets $A_1,\dots,A_s$. Define
\begin{eqnarray}\label{Ex:mixed-vol-val}
\phi(K)=V(K[n-s],A_1,\dots,A_s)
\end{eqnarray}
where $K[n-s]$ means that $K$ is taken $n-s$ times. Then $\phi$ is a
translation invariant continuous valuation. This easily follows from
Example \ref{Ex:valuations}(3).

\subsection{Hadwiger's theorem.}\label{Ss:hadwiger-thm}
One of the most famous and classical results of the valuations
theory is Hadwiger's classification of isometry invariant continuous
convex valuations on the Euclidean space $\RR^n$. To formulate it,
let us denote by $V_i$ the $i$-th intrinsic volume, which by
definition is
$$V_i(K)=c_{n,i}V(K[i],D[n-i])$$
where $c_{n,i}$ is an explicitly written constant which is just a
standard normalization (see \cite{schneider-book}). In particular
$V_0=\chi$ is the Euler characteristic, $V_n=vol$ is the Lebesgue
measure normalized so that  the volume of the unit cube is equal to
1. Clearly $V_i\in Val_i$ is an $O(n)$-invariant valuation.
\begin{theorem}[Hadwiger's classification \cite{hadwiger-book}]\label{T:hadwiger}
Any $SO(n)$-invariant translation invariant continuous convex
valuation is a linear combination of $V_0,V_1,\dots,V_n$. (In
particular it is $O(n)$-invariant.)
\end{theorem}
In 1995 Klain \cite{klain-hadwiger} has obtained a simplified proof
of this deep result as an easy consequence of his classification of
simple even valuations discussed below in Section
\ref{Ss:klain-schneider}. Hadwiger's theorem turned out to be very
useful in integral geometry of the Euclidean space. This will be
discussed in more detail in J. Fu's lectures. We also refer to the
book \cite{klain-rota}.

\subsection{Irreducibility theorem.}\label{Ss:irreducibility}
One of the basic questions of the valuations theory is to describe
valuations with given properties. Hadwiger's theorem is one example
of such a result of great importance. In recent years there were
obtained a number of classification results of various classes of
valuations. The case of continuous translation invariant valuations
will be discussed in detail in these lectures below and in lectures
by Fu.

The question is whether it is possible to give a reasonable
description of all translation invariant continuous convex
valuations. In 1980 P. McMullen \cite{mcmullen-80} has formulated a
more precise conjecture which says that linear combinations of mixed
volumes (as in (\ref{Ex:mixed-vol-val})) are dense in $Val$. This
conjecture was proved in positive by the author
\cite{alesker-gafa-01} in a stronger form which later on turned out
to be important in further developments and applications.

To describe the result let us make a few more remarks. We say that a
valuation $\phi$ is even (resp. odd) if $\phi(-K)=\phi(K)$ (resp.
$\phi(-K)=-\phi(K)$) for any $K\in \ck(V)$. The subspace of even
(resp. odd) $i$-homogeneous valuations will be denoted by $Val_i^+$
(resp. $Val_i^-$). Clearly
\begin{eqnarray}\label{E:parity-decomposition}
Val_i=Val_i^+\oplus Val_i^-.
\end{eqnarray}

Next observe that the group $GL(V)$ of all invertible linear
transformations acts linearly on $Val$:
$$(g\phi)(K)=\phi(g^{-1}K).$$
\begin{theorem}[Irreducibility theorem
\cite{alesker-gafa-01}]\label{T:irreducibility} For each $i$, the
spaces $Val_i^+,\, Val_i^-$ are irreducible representations of
$GL(V)$, i.e. they do not have proper invariant closed subspaces.
\end{theorem}

\begin{remark}
By Remark \ref{R:classification-extremal-cases} $Val_0^+=Val_0$,
$Val^+_n=Val_n$ are one dimensional. But for $1\leq i\leq n-1$ the
spaces $Val_i^{\pm}$ are infinite dimensional. $Val_{n-1}$ was
explicitly described by McMullen \cite{mcmullen-80}; we state his
result in Section \ref{Ss:klain-schneider} below.
\end{remark}

Theorem \ref{T:irreducibility} immediately implies McMullen's
conjecture. Indeed it is easy to see that the closure of the linear
span of mixed volumes is a $GL(V)$-invariant subspace, and its
intersection with any $Val_i^\pm$ is non-zero. Hence by the
irreducibility theorem any such intersection should be equal to the
whole space $Val_i^\pm$.

The irreducibility theorem will be used in these lectures several
times. The proof of this result uses a number of deep results from
the valuations theory in combination with representation theoretical
techniques. A particularly important such result of high independent
interest is the Klain-Schneider classification of {\itshape simple}
translation invariant continuous convex valuations; it is discussed
in the next section.

\subsection{Klain-Schneider characterization of simple
valuations.}\label{Ss:klain-schneider}
\begin{definition}\label{D:simple-valuations}
A convex valuation $\phi\in Val$ is called simple if $\phi(K)=0$ for
any $K\in \ck(V)$ with $\dim K<n:=\dim V$.
\end{definition}
\begin{theorem}\label{T:kl-sch}
(i) [Klain \cite{klain-hadwiger}] Any simple even valuation from
$Val$ is proportional to the Lebesgue measure.

(ii)[Schneider \cite{schneider-simple}] Any simple odd valuation
from $Val$ is $(n-1)$-homogeneous.
\end{theorem}

Clearly any simple valuation is the sum of a simple even and a
simple odd valuations. Hence in order to complete the description of
simple valuations it remains to classify simple $(n-1)$-homogeneous
valuations. Fortunately McMullen \cite{mcmullen-80} has previously
described $Val_{n-1}$ very explicitly. His result was used in
Schneider's proof, and it is worthwhile to state it explicitly as it
has independent interest.

First let us recall the definition of the area measure
$S_{n-1}(K,\cdot)$ of a convex compact set $K$. Though it is not
strictly necessary, it is convenient and common to fix a Euclidean
metric on $V$. After this choice, $S_{n-1}(K,\cdot)$ is a measure on
the unit sphere $S^{n-1}$ defined as follows. First let us assume
that $K$ is a polytope. For any $(n-1)$-face $F$ let us denote by
$n_F$ the unit outer normal at $F$. Then by definition
$$S_{n-1}(K,\cdot)=\sum_F vol_{n-1}(F)\delta_{n_F},$$
where the sum runs over all $(n-1)$-faces of $K$, and $\delta_{n_F}$
denotes the delta-measure supported at $n_F$. Then the claim is that
the area measure extends uniquely by weak continuity to the class of
all convex compact sets: if  $K_N\to K$ in the Hausdorff metric then
$S_{n-1}(K_N\cdot)\to S_{n-1}(K,\cdot)$ weakly in the sense of
measures (see \cite{schneider-book}, \S 4.2).

\begin{theorem}[McMullen, \cite{mcmullen-80}]\label{mcmullen-classif}
Let $\phi\in Val_{n-1}$, $n=\dim V$. Then there exists a continuous
function $f\colon S^{n-1}\to \CC$ such that
$$\phi(K)=\int_{S^{n-1}}f(x)dS_{n-1}(K,x).$$
Conversely, any expression of this form with a continuous $f$ is a
valuation from $Val_{n-1}$.

Moreover two continuous functions $f$ and $g$ define the same
valuation if and only if the difference $f-g$ is a restriction of a
linear functional on $V$ to the unit sphere.
\end{theorem}

Now we can summarize the classification of simple valuations.
\begin{theorem}[Klain-Schneider]\label{T:klain-schneider}
Simple translation invariant continuous valuations are precisely of
the form
$$K\mapsto c\cdot vol_n(K)+\int_{S^{n-1}}f(x)dS_{n-1}(K,x),$$
where $f\colon S^{n-1}\to \CC$ is an odd continuous function, and
$c$ is a constant. Moreover the constant $c$ is defined uniquely,
while $f$ is defined up to a linear functional.
\end{theorem}

\subsection{Smooth translation invariant
valuations.}\label{Ss:smooth-val} We are going to describe an
important subspace of $Val$ of {\itshape smooth} valuations. They
form a dense subspace is $Val$ and carry a number of extra
structures (e.g. product, convolution, Fourier transform) which do
not extend by continuity to the whole space $Val$ of continuous
valuations. Moreover main examples relevant to integral geometry are
in fact smooth valuations.

\begin{definition}\label{D:smooth-transl}
A valuation $\phi\in Val(V)$ is called smooth if the Banach space
valued map $GL(V)\to Val(V)$ given by $g\mapsto g(\phi)$ is
infinitely differentiable.
\end{definition}

From a very general and elementary representation theoretical
reasoning, the subset of smooth valuations, denoted by
$Val^{sm}(V)$, is a linear dense subspace of $Val(V)$ invariant
under the natural action of $GL(V)$. Also $Val^{sm}(V)$ carries a
linear topology which is stronger than that induced from $Val(V)$,
and with respect to which it is a Fr\'echet space. This is called
often the Garding topology, and tacitly we will always assume that
$Val^{sm}(V)$ is equipped with it. Of course, $Val^{sm}$ also
satisfies McMullen's decomposition and the irreducibility theorem.

For future applications to integral geometry, the following result
will be important.

\begin{proposition}[\cite{alesker-gafa-04}]\label{P:smooth-invariant-val}
Let $G$ be a compact subgroup of the orthogonal group of a Euclidean
space $V$. Assume that $G$ acts transitively on the unit sphere of
$V$. Then any $G$-invariant valuation from $Val(V)$ is smooth.
\end{proposition}

\subsection{Product on smooth translation invariant
valuations and Poincar\'e duality.}\label{Ss:product} In this
section we discuss the product on translation invariant smooth
valuations introduced in \cite{alesker-gafa-04}. This structure
turned out to be intimately related to integral geometric formulas
discussed in detail in J. Fu's lectures.

First we will introduce exterior product on valuations.

\begin{theorem}[\cite{alesker-gafa-04}]\label{T:ext-product}
Let $V$ and $W$ be finite dimensional real vector spaces. There
exists a continuous bilinear map, called exterior product,
$$Val^{sm}(V)\times Val^{sm}(W)\to Val(V\times W)$$
which is uniquely characterized by the following property: Fix $A\in
\ck(V),\, B\in \ck(W)$. Let $vol_V,\, vol_W$ be Lebesgue measures on
$V,W$ respectively. Let $\phi(K)=vol_V(K+A),\, \psi(K)=vol_W(K+B)$.
Then their exterior product, denoted by $\phi\boxtimes \psi$, is
$$(\phi\boxtimes\psi)(K)=vol_{V\times W}(K+(A\times B)) \mbox{ for any } K\in \ck(V\times W),$$
where $vol_{V\times W}$ is the product measure of $vol_V$ and
$vol_W$.

%(2) The exterior product is associative:
%$$\phi_1\boxtimes(\phi_2\boxtimes\phi_3)=(\phi_1\boxtimes\phi_2)\boxtimes\phi_3.$$
\end{theorem}
Notice that the uniqueness in this theorem follows immediately from
the (proved) McMullen's conjecture since linear combinations of
valuations of the form $vol(\bullet +A)$ are dense in valuations.

Let us emphasize that the exterior product is defined on smooth
valuations, but it takes values not in smooth but just continuous
valuations. Usually the exterior product is not smooth. Let us give
some examples.
\begin{example}
1) The exterior product of Lebesgue measures in the sense of
valuations coincides obviously with their measure theoretical
product.

2) The exterior product of Euler characteristics is the Euler
characteristic on $V\times W$.

3) Let $vol_V$ be a Lebesgue measure on $V$, and $\chi_W$ be the
Euler characteristic on $W$. Then the exterior product
$vol_V\boxtimes \chi_W$ is the volume of the projection to $V$:
$$(vol_V\boxtimes \chi_W)(K)=vol_V(pr_V(K)) \mbox{ for any } K\in \ck(V\times W),$$
where $pr_V\colon V\times W\to V$ is the natural projection. Observe
that this valuation is not smooth (in contrast to the first two
examples.)
\end{example}

The first non-trivial point in Theorem \ref{T:ext-product} is that
the exterior product is well defined, the second one is continuity.
We do not give here any proof. However let us give an incomplete,
but instructive, explanation why the exterior product is well
defined. There are of course many different ways to write a
valuation as a linear combinations of $vol(\bullet +A)$. Let us see
that the exterior product of finite linear combinations of such
expressions is independent of the particular choice of a linear
combination. Since the situation is symmetric with respect to both
valuations, we may assume that $\phi(\bullet)=\sum_ic_i\cdot
vol_V(\bullet +A_i)$ and $\psi (\bullet)=vol_W(\bullet +B)$. Then
using the Fubini theorem and the equality $A_i\times B=(A_i\times
0)+(0\times B)$ we get
\begin{eqnarray*}
(\phi\boxtimes\psi)(K)=\sum_ic_i\cdot vol_{V\times W}(K+(A_i\times
B))=\\
\sum_i c_i\cdot \int_{w\in W}vol_{V}\left[ \left\{(K+(0\times
B))\cap
(V\times \{w\})\right\}+A_i\right] dvol_W(w)=\\
\int_{w\in W} \phi\left[(K+(0\times B))\cap (V\times \{w\})\right]
dvol_W(w).
\end{eqnarray*}
Clearly the last expression is independent of the form of
presentation of $\phi$.

\hfill

Now let us define the product on $Val^{sm}$. Let us denote by
$$\Delta\colon V\to V\times V$$
the diagonal imbedding. The product of $\phi,\psi\in Val^{sm}(V)$ is
defined by
$$(\phi\cdot \psi):=(\phi\boxtimes \psi)(\Delta(K)).$$
It turns out that the product of smooth valuations is again smooth.
\begin{theorem}[\cite{alesker-gafa-04}]
The product of smooth valuations $Val^{sm}(V)\times Val^{sm}(V)\to
Val^{sm}(V)$ is continuous (in the Garding topology), associative,
commutative, and distributive. Then $Val^{sm}(V)$ becomes an algebra
over $\CC$ with unit which is the Euler characteristic. Moreover the
product respects the degree of homogeneity: $$Val^{sm}_i\cdot
Val_j^{sm}\subset Val_{i+j}^{sm}.$$
\end{theorem}

\begin{example}\label{Ex:product-intrinsic}
The product of intrinsic volumes $V_i\cdot V_j$ with $i+j\leq n$ is
a non-zero multiple of $V_{i+j}$: by the Hadwiger theorem it is
clear that the product should be proportional to $V_{i+j}$, the
constant of proportionality can be computed explicitly.
\end{example}

An interesting property of the above product is a version of the
Poincar\'e duality.
\begin{theorem}[\cite{alesker-gafa-04}]\label{T:poincare-dual-finite}
For any $i=0,1,\dots,n=\dim V$ the bilinear map
$$Val^{sm}_i(V)\times Val^{sm}_{n-i}(V)\to Val_n^{sm}(V)$$
is a perfect pairing, namely for any non-zero valuation $\phi\in
Val^{sm}_i(V)$ there exists $\psi\in Val^{sm}_{n-i}(V)$ such that
$\phi\cdot\psi\ne 0$.
\end{theorem}
This result follows easily from the Irreducibility Theorem
\ref{T:irreducibility}. Indeed it suffices to prove the statement
for valuations of fixed parity $\eps=\pm 1$. Then the kernel of the
above pairing in $Val^{\eps,sm}_i(V)$ is a $GL(V)$-invariant closed
subspace. Hence it must be either zero or everything. But it cannot
be everything since then for any  valuation $\psi\in
Val^{sm}_{n-i}(V)$ one would have $\psi\cdot Val_i^{\eps,sm}(V)=0$.
But this is not the case as can be easily proved by constructing an
explicit example. (Say in the even case, the product of the
intrinsic volumes $V_i\cdot V_{n-i}$ is a non-zero multiple of
Lebesgue measure.)

Thus $Val^{sm}(V)$ is a graded algebra satisfying the Poincar\'e
duality. In Section \ref{Ss:hard-lefschetz} we will see also that it
satisfies two versions of the hard Lefschetz theorem.

\subsection{Pull-back and push-forward of translation invariant
valuations.}\label{Ss:pull-push-translation} In this section we
describe the operations of pull-back and push-forward on translation
invariant valuations under linear maps.

Let $f\colon V\to W$ be a linear map. We define
\cite{alesker-fourier} the continuous linear map, called pull-back,
$$f^*\colon Val(W)\to Val(V)$$
defined as usual by $(f^*\phi)(K)=\phi(f(K))$. It is easy to see
that $f^*\phi$ is indeed a continuous translation invariant convex
valuation. The following result is evident.
\begin{proposition}
(i) $f^*$ preserves the degree of homogeneity and parity.

(ii) $(f\circ g)^*=g^*\circ f^*$.
\end{proposition}

Notice that the product on valuations can be expressed via the
exterior product and the pull-back by
$$\phi\cdot \psi=\Delta^*(\phi\boxtimes\psi),$$
where $\Delta$ is the diagonal imbedding.

\hfill

A somewhat less obvious operation is push-forward $f_*$ introduced
in \cite{alesker-fourier}. In some non-precise sense $f_*$ is dual
to $f^*$. In these notes it will be used only to give an alternative
description of the convolution on valuations in Section
\ref{Ss:convolution} and to clarify some properties of the Fourier
type transform on valuations in Section \ref{Ss:fourier}; the reader
not interested in these subjects may skip the rest of this section.

Canonically the push-forward map acts not between spaces of
valuations, but between their tensor product (twist) by an
appropriate one dimensional space of Lebesgue measures. To be more
precise let us denote by $D(V^*)$ the one dimensional space of
($\CC$-valued) Lebesgue measures on $V^*$. Then $f_*$ is a linear
continuous map
$$f_*\colon Val(V)\otimes D(V^*)\to Val(W)\otimes D(W^*).$$
In order to define this map, we will split its construction to the
cases of $f$ being injective, surjective, and a general linear map.

\underline{Case 1. Let $f$ be injective.} Thus we may assume that
$V\subset W$. In order to simplify the notation we choose a
splitting $W=V\oplus L$ and we fix Lebesgue measures on $V$ and $L$.
Then on $W$ we have the product measure. These choices induce
isomorphisms $Dens(V^*)\simeq \CC $, $Dens(W^*)\simeq \CC$. We leave
for a reader to check that the construction of $f_*$ is independent
of these choices.

Let $\phi\in Val(V)$. Define
$$(f_*\phi)(K)=\int_{l\in L}\phi(K\cap (l+V))dl.$$
It is easy to see that $f_*\colon Val(V)\to Val(W)$ is a continuous
linear map.

\underline{Case 2. Let $f$ be surjective.} Again it will be
convenient to assume that $f$ is just a projection to a subspace,
and fix a splitting $V=W\oplus M$. Again fix Lebesgue measures on
$W,M$, and hence on $V$. Let us also fix a set $S\in \ck(M)$ of unit
measure. Set $m:=\dim M$. Then define
\begin{eqnarray*}
(f_*\phi)(K)=\frac{1}{m!}\frac{d^m}{d\eps^m}\phi(K+\eps \cdot
S)\big|_{\eps=0}.
\end{eqnarray*}
Recall that by McMullen's theorem $\phi(K+\eps\cdot S)$ is a
polynomial in $\eps \geq 0$. Moreover its degree is at most $m$:
Indeed when $K$ is fixed this expression is a translation invariant
continuous valuation with respect to $S\in \ck(M)$. The coefficient
of $\eps^m$ is an $m$-homogeneous valuation with respect to
$S\subset M$, and hence by Hadwiger's theorem (see Remark
\ref{R:classification-extremal-cases}) it must be proportional to
the Lebesgue measure on $M$ with a constant depending on $K$. By our
definition, this coefficient is exactly $(f_*\phi)(K)$, in
particular it does not depend on $S$. In fact it does not depend
also on choice of Lebesgue measures and the splitting.

\underline{Case 3. Let $f$ be a general linear map.} Let us choose a
factorization $f=g\circ h$ where $h\colon V\to Z$ is surjective, and
$g\colon Z\to W$ is injective. Then define $f_*:=g_*\circ h_*$. One
can show that $f_*$ is independent of the choice of such a
factorization.

\begin{proposition}[\cite{alesker-fourier},
Section 3.2]\label{P:push-forward-transl} (i) $f_*\colon
Val(V)\otimes D(V^*)\to Val(W)\otimes D(W^*)$ is a continuous linear
operator.

(ii) $(f\circ g)_*=f_*\circ g_*$.

(iii) $f_*\left(Val_i(V)\otimes D(V^*)\right)\subset Val_{i-\dim
V+\dim W}(W)\otimes D(W^*)$.
\end{proposition}

\subsection{Convolution.}\label{Ss:convolution} In this section we
describe another interesting operation on valuations: a convolution
introduced by Bernig and Fu \cite{bernig-fu-convolution}. This is a
continuous product on $Val^{sm}\otimes D(V^*)$. Let us fix for
simplicity of the notation a Lebesgue measure $vol$ on $V$; it
induces a Lebesgue measure on $V^*$. With these identifications,
convolution is going to be defined on $Val^{sm}(V)$ (without the
twist by $D(V^*)$).
\begin{theorem}[\cite{bernig-fu-convolution}]
There exists a unique continuous bi-linear map, called convolution,
$$\ast\colon Val^{sm}(V)\times  Val^{sm}(V)\to Val^{sm}(V)$$ such that
$$vol(\bullet +A)\ast vol(\bullet +B)=vol(\bullet +A+B).$$
This product makes $Val^{sm}(V)$ a commutative associative algebra
with the unit element $vol$. Moreover $Val^{sm}_i\ast
Val^{sm}_j\subset Val^{sm}_{i+j-n}$.
\end{theorem}

The above result characterizes the convolution uniquely, and allows
to compute it in some examples. We can give however one more
description of it using the previously introduced operations. Let
$a\colon V\times V\to V$ be the addition map, namely $a(x,y)=x+y$.
Then by \cite{alesker-fourier}, Proposition 3.3.2, one has
\begin{eqnarray*}
\phi\ast\psi=a_*(\phi\boxtimes \psi).
\end{eqnarray*}

The product and convolution will be transformed one to the other in
Section \ref{Ss:fourier} by another useful operation, the Fourier
type transform.

\subsection{Hard Lefschetz type theorems.}\label{Ss:hard-lefschetz}
The product and the convolution on valuations satisfy another
non-trivial property analogous to the hard Lefschetz theorem from
algebraic geometry \cite{griffiths-harris}. Let us fix on $V$ a
Euclidean metric. Consider the operator
$$L\colon Val^{sm}_\ast\to Val^{sm}_{\ast+1}$$
given by $L\phi:=\phi\cdot V_1$ where $V_1$ is the first intrinsic
volume as in Section \ref{Ss:hadwiger-thm}. Consider also another
operator
$$\Lam\colon Val^{sm}_\ast\to Val^{sm}_{\ast-1}$$ defined by
$(\Lam \phi)(K)=\frac{d}{d\eps}\phi(K+\eps \cdot D)\big|_{\eps=0}$
where $D$ is the unit ball (here we use again McMullen's theorem
that $\phi(K+\eps\cdot D)$ is a polynomial). Notice that up to a
normalizing constant, the operator $\Lam$ is equal to the
convolution with $V_{n-1}$, as was observed by Bernig and Fu
\cite{bernig-fu-convolution}.
\begin{theorem}\label{T:hard-lefschetz}
(i) Let $0\leq i<\frac{n}{2}$. Then $L^{n-2i}\colon Val^{sm}_i\to
Val^{sm}_{n-i}$ is an isomorphism.

(ii) Let $\frac{n}{2}<i\leq n$. Then $\Lam^{2i-n}\colon
Val^{sm}_i\to Val^{sm}_{n-i}$ is an isomorphism.
\end{theorem}
Several authors have contributed to the proof of this theorem. First
the author proved (i) and (ii) in the even case
\cite{alesker-gafa-sem-04}, \cite{alesker-jdg-03} using the previous
joint work with Bernstein \cite{alesker-bernstein} and integral
geometry on Grassmannians (Radon and cosine transforms). Then Bernig
and Br\"ocker \cite{bernig-brocker} proved part (ii) in the odd case
using a very different method: the Laplacian acting on differential
forms on the sphere bundle and some results from complex geometry
(K\"ahler identities). Next Bernig and Fu have shown
\cite{bernig-fu-convolution} that in the even case, both versions of
the hard Lefschetz theorem are in fact equivalent via the Fourier
transform (which was at that time defined only for even valuations).
Finally the author extended in \cite{alesker-fourier} the Fourier
transform to odd valuations and derived the version (i) of the hard
Lefschetz theorem in the odd case from version (ii).

\subsection{A Fourier type transform on translation invariant convex
valuations.}\label{Ss:fourier} A Fourier type transform on
translation invariant smooth valuations is another useful operation.
First it was introduced in \cite{alesker-jdg-03} (under a different
name of duality) for even valuations and was applied there to
Hermitian integral geometry in order to construct a new basis in the
space of $U(n)$-invariant valuations on $\CC^n$. In the odd case it
was constructed in \cite{alesker-fourier}. Some recent applications
and non-trivial computations of the Fourier transform in Hermitian
integral geometry are due to Bernig and Fu \cite{bernig-fu-annals}.

In this section we will describe the general properties of the
Fourier transform and its relation to the product and convolution
described above. We will present a construction of the Fourier
transform in the even case only. The construction in the odd case is
more technical, and will not be presented here. Notice that the even
case will suffice for a reader interested mostly in applications to
integral geometry of affine spaces, since by a result of Bernig
\cite{bernig-spin7-g2} any $G$-invariant valuation from $Val$ must
be even provided $G$ is a compact subgroup of the orthogonal group
acting transitively on the unit sphere.

The main general properties of the Fourier transform are summarized
in  the following theorem.
\begin{theorem}[\cite{alesker-fourier}]\label{T:fourier}
There exists an isomorphism of linear topological spaces
$$\FF\colon Val^{sm}(V)\to Val^{sm}(V^*)\otimes D(V)$$ which
satisfies the following properties:

1) $\FF$ commutes with the natural action of the group $GL(V)$ on
both spaces;

2) $\FF$ is an isomorphism of algebras when the source is equipped
with the product and the target with the convolution.

3) The Fourier transform satisfies a Plancherel type inversion
formula explained below.
\end{theorem}

In order to describe the Plancherel type formula and present a more
explicit description of the Fourier transform it will be convenient
(but not necessary) to fix a Euclidean metric on $V$. This will
induce identifications $V^*\simeq V$ and $D(V^*)\simeq \CC$. With
these identifications $\FF\colon Val^{sm}(V)\to Val^{sm}(V)$;
actually it changes the degree of homogeneity:
$$\FF\colon Val_i^{sm}\tilde\to Val_{n-i}^{sm}.$$
The Plancherel type formula says, under these identifications, that
$(\FF^2\phi)(K)=\phi(-K)$.

Here are a few simple examples: $\FF(vol)=\chi$; $\FF(\chi)=vol$;
$\FF(V_i)=c_{n,i}V_{n-i}$ where $c_{n,i}>0$ is a normalizing
constant which can be computed explicitly. (Notice that the last
fact, except for the positivity of $c_{n,i}$, is an immediate
consequence of the fact that $\FF$ commutes with the action of
$O(n)$ and Hadwiger's theorem.)

\hfill

%**********************************

The Fourier transform on a 2-dimensional space has an explicit
description which we are going to describe now. We will work for
simplicity in $\RR^2$ with the standard Euclidean metric and
standard orientation. With the identifications induced by the metric
as above $\FF\colon Val^{sm}(\RR^2)\to Val^{sm}(\RR^2)$. It remains
to describe $\FF$ on 1-homogeneous valuations. Every such smooth
valuation $\phi$ can be written uniquely in the form
$$\phi(K)=\int_{S^1}h(\ome) dS_1(K,\ome)$$
where $h\colon S^1\to \CC$ is a smooth function which is orthogonal
on $S^1$ to the two dimensional space of linear functionals. Let us
decompose $h$ to the even and odd parts:
$$h=h_++h_-.$$
Let us decompose further the odd part $h_-$ to "holomorphic" and
"anti-holomorphic" parts
$$h_-=h_-^{hol}+h_-^{anti}$$
as follows. Let us decompose $h_-$ to the usual Fourier series on
the circle $S^1$:
$$h_-(\ome)=\sum_{k} \hat h_-(k)e^{ik\ome}.$$
Then by definition
\begin{eqnarray*}
h_-^{hol}(\ome):=\sum_{k>0} \hat h_-(k)e^{ik\ome},\\
h_-^{anti}(\ome):=\sum_{k<0} \hat h_-(k)e^{ik\ome}.
\end{eqnarray*}
Then the Fourier transform of the valuation $\phi$ is equal to
$$(\FF\phi)(K)=\int_{S^1}(h_+(J\ome)+h_-^{hol}(J\ome))dS_1(K,\ome)-\int_{S^1}h_-^{anti}(J\ome)dS_1(K,\ome)$$
where $J$ is the rotation of $\RR^2$ by $\frac{\pi}{2}$
counterclockwise. (Notice the minus sign before the second
integral.) Observe that $\FF$ preserves the class of real valued
even valuations, but for odd real valued valuations this is not
true. This phenomenon also holds in higher dimensions.

%*****************************

\hfill

Let us consider even valuations in arbitrary dimension. We again fix
a Euclidean metric on $V$. A useful tool in studying even valuations
is an imbedding of $Val^+_i(V)$ to the space of continuous functions
on the Grassmannian $Gr_i(V)$ of $i$-dimensional subspaces of $V$.
It was constructed by Klain \cite{klain-even} as an easy consequence
of his classification of simple even valuation (Theorem
\ref{T:kl-sch}(a)). Define the map
$$Kl_i\colon Val^+_i(V)\to C(Gr_i(V))$$
as follows. Let $\phi\in Val^+_i(V)$. For any $E\in Gr_i(V)$ the
restriction of $\phi$ to $E$ is a valuation of maximal degree of
homogeneity. Hence by a result of Hadwiger it must be proportional
to Lebesgue measure $vol_E$ induced by the Euclidean metric. Thus by
definition
$$\phi|_E= \left(Kl_i(\phi)\right)(E)\cdot vol_E.$$
Clearly $Kl_i(\phi)$ is a continuous function and $Kl_i$ is a
continuous linear $O(n)$-equivariant linear map. The non-trivial
fact is that $Kl_i$ is {\itshape injective}. For we observe that if
$\phi\in Ker (Kl_i)$ then the restriction of $\phi$ to any
$i+1$-dimensional subspace is a simple, even, $i$-homogeneous
valuation. Hence it vanishes by Klain's theorem. Proceeding by
induction, one sees that $\phi=0$.

Next it is not hard to see that smooth valuations are mapped under
$Kl_i$ to infinitely smooth functions on $Gr_i(V)$. Let us define
the Fourier transform $\FF\colon Val_i^{+,sm}\to Val_{n-i}^{+,sm}$
by the following property: for any subspace $E\in Gr_{n-i}(V)$,
$$(Kl_{n-i}(\FF\phi))(E)=(Kl_i(\phi))(E^\perp),$$
where as usual $E^\perp$ denotes the orthogonal complement. This
condition characterizes $\FF$ uniquely in the even case. The
non-trivial point is the existence of $\FF$ with this property. The
problem is that the Klain imbedding $Kl_i\colon Val_i^{+,sm}(V)\inj
C^\infty(Gr_i(V))$ is not onto (for $i\ne 1,n-1$). The main point is
to show that this image is invariant under taking the orthogonal
complement. It was shown by Bernstein and the author
\cite{alesker-bernstein} that the image of $Kl_i$ coincides with the
image of the so called cosine transform on Grassmannians;  this step
used also the irreducibility theorem. From the definition of the
cosine transform (which we do not reproduce here) it is easy to see
that its image is invariant under taking the orthogonal complement.

\hfill

Let us add a couple of words on the odd case. There is a version of
Klain's imbedding for odd valuations though  it is more complicated:
$Val_i^{-,sm}(V)$ is realized as a quotient of a subspace of
functions on a manifold of partial flags (here instead of Klain's
characterization of simple even valuations one has to use
Schneider's version for odd valuations - Theorem \ref{T:kl-sch}(b)).
We call it Schneider's imbedding. However there is no direct
analogue of "the cosine transform" description of the image of it.
More delicate analysis is required; it is based (besides the
irreducibility theorem) on a more detailed study of the action of
$GL(n,\RR)$ on valuations and on functions (or, rather, sections of
an equivariant line bundle) on partial flags. This requires more
tools from infinite dimensional representation theory of the group
$GL(n,\RR)$. We refer for the details to \cite{alesker-fourier}.

\hfill

Another important property of the Fourier transform is that it
intertwines the pull-back and push-forward on valuations. We will
formulate this property in a non-rigorous way due to various
technicalities making the formal statement heavier (see
\cite{alesker-integral}). Let $f\colon V\to W$ be a linear map. Let
$f^\vee\colon W^*\to V^*$ be the dual map between the dual spaces.
Then the claim is that one should have
\begin{eqnarray}\label{E:diagram-fourier}
\FF_V\circ f^*=(f^\vee)_*\circ \FF_W,
\end{eqnarray} where $f^*$ is
the pull-back under $f$, $(f^\vee)_*$ is the push-forward under
$f^\vee$, and $\FF_V$ and $\FF_W$ are the Fourier transforms on $V$
and $W$ respectively. Notice that the equality
(\ref{E:diagram-fourier}) formally is ill-defined if $f$ is not an
isomorphism. This is due to the fact that $\FF$ is formally defined
only on the class of smooth valuations, while $f^*$ and $(f^\vee)_*$
do not preserve this class. Nevertheless morally this equality
should be true, but technically one should be more accurate.

\hfill

Moreover one expects that in some sense the Fourier transform should
commute with the exterior product:
$$\FF(\phi\boxtimes\psi)=\FF\phi\boxtimes \FF\psi.$$
The difficulty here is that the exterior product of smooth
valuations is usually not smooth.

\hfill

As the last remark let us mention that it would be desirable to have
a more direct construction of the Fourier transform. For example we
still do not know how to describe it in terms of the construction of
valuations using integration with respect to the normal cycle
discussed below in Section \ref{Sss:normal-cycle}.

\subsection{General constructions of translation invariant
convex valuations.}\label{Ss:constructions} So far the only
construction of valuations we have seen is the mixed volume. In this
section we review some other general constructions of translation
invariant continuous convex valuations. In Section
\ref{Sss:integral-geometry} we describe briefly an array of examples
coming from integral geometry; more complete treatment will be given
in Fu's lectures. In Section \ref{Sss:normal-cycle} we describe
another very general and useful construction via integration over
the normal cycle of a set; this construction will be generalized
appropriately to the context of valuations on manifolds in Section
\ref{S:valuations-on-manifolds}. There is yet another construction
of valuations based on complex and quaternionic pluripotential
theory. It is somewhat more specialized and will not be discussed
here; we refer to \cite{alesker-adv-05}, \cite{alesker-spin9}, and
the survey \cite{alesker-survey-alg-analys}.

\subsubsection{Integral geometry.}\label{Sss:integral-geometry}
Let us give a few basic examples which arise naturally in (Crofton
style) integral geometry. The classical reference to this type of
integral geometry is Santal\'o's book \cite{santalo-book}. For
further discussion of this type of integral geometry and its
relations to the valuations theory we refer to Fu's lectures, the
book \cite{klain-rota}, and the articles \cite{alesker-jdg-03},
\cite{bernig-sun}, \cite{bernig-spin7-g2},
\cite{bernig-quaternions}, \cite{bernig-fu-annals},
\cite{fu-unitary} (these recent results are surveyed by  Bernig
\cite{bernig-survey}).

\hfill

Let $V=\RR^n$ be the standard Euclidean space. Let $Gr_{k,n}$ denote
the Grassmannian of all linear $k$-dimensional subspaces of it, and
let $\bar Gr_{k,n}$ denote the Grassmannian of affine
$k$-dimensional subspaces.
%Let $O(n)$ denote the groups of all linear isometries of
%$\RR^n$, and $\bar O(n)$ denote the group of all affine isometries.
It is not hard to check that the following expressions are
continuous valuations invariant with respect to all isometries of
$\RR^n$:
\begin{eqnarray}\label{E:int-euclid1}
\phi(K)=\int_{E\in Gr_{k,n}}V_i(pr_E(K))dE,\\\label{E:int-euclid2}
\psi(K)=\int_{E\in \bar Gr_{k,n}} V_i(K\cap E)dE,
\end{eqnarray}
where $dE$ denotes in both formulas a Haar measure on the
corresponding Grassmannian, and $pr_E\colon \RR^n\to E$ denotes the
orthogonal projection. These expressions have been studied quite
extensively in the classical integral geometry; they can be computed
as integrals of certain expressions of the principal curvatures of
the boundary $\pt K$, at least under appropriate assumptions on
smoothness of it, see e.g. \cite{santalo-book}, \cite{klain-rota}.
Notice that Hadwiger's theorem implies that these valuations can be
written as linear combinations of intrinsic volumes
$V_0,V_1,\dots,V_n$; the coefficients can be computed explicitly.

\hfill

Let us present analogous expressions from the Hermitian integral
geometry of $\CC^n$. Despite the obvious similarity to the Euclidean
case, these expressions have been studied in depth only quite
recently \cite{alesker-jdg-03}, \cite{bernig-fu-annals},
\cite{fu-unitary}.
%We denote by $U(n)$ the group of unitary transformations, and by $\bar U(n)$
%the group $U(n)\ltimes \CC^n$ of holomorphic isometries.
Let $\grc_{k,n}$ denote the Grassmannian of complex linear
$k$-dimensional subspaces of $\CC^n$, and $\bar \grc_{k,n}$ the
Grassmannian of complex affine $k$-dimensional subspaces. Let us
define in analogy to (\ref{E:int-euclid1})-(\ref{E:int-euclid2})
\begin{eqnarray}\label{E:int-hermit1}
\phi(K)=\int_{E\in \grc_{k,n}}V_i(pr_E(K))dE,\\\label{E:int-hermit2}
\psi(K)=\int_{E\in \bar\grc_{k,n}}V_i(K\cap E)dE,
\end{eqnarray}
where $dE$ again denotes a Haar measure on the appropriate complex
Grassmannian. It was shown in \cite{alesker-jdg-03} that from
valuations of the form (\ref{E:int-hermit1}) (or alternatively,
(\ref{E:int-hermit2})) one can choose a basis of unitarily invariant
valuations in $Val(\CC^n)$. Moreover in the same paper it was shown
that the Fourier transform of a valuation of the form
(\ref{E:int-hermit1}) has the form (\ref{E:int-hermit2}) with
appropriately chosen $i,k$, and vise versa. Some different bases in
unitarily invariant valuations have been constructed by Bernig-Fu
\cite{bernig-fu-annals} where they also computed several integral
geometric formulas in $\CC^n$, in particular the principal kinematic
formula.

\subsubsection{Normal cycle.}\label{Sss:normal-cycle} In this
section we remind the notion of the normal cycle of a convex set and
present a construction of translation invariant smooth valuations on
convex sets in terms of it. In fact we will see that all such
valuations can be obtained using this construction. One of the
important aspects of this construction  is that it generalizes to a
broader context of valuations on manifolds to be discussed in
Section \ref{S:valuations-on-manifolds}.

\hfill

In this section we will fix again a Euclidean metric and an
orientation on a vector space $V$, $\dim V=n$, for convenience of a
geometrically oriented reader. However this metric is not necessary,
and in Section \ref{Ss:definitions-manifolds} we describe an
extension of the construction of normal cycle to any smooth manifold
without any additional structure (not for convex sets of course, but
for compact submanifolds with corners).

Let $K\in \ck(V)$ be a convex compact subset of $V$. For any point
$x\in K$ let us define the normal cone of $K$ at $x$ as a subset of
the unit sphere $S^{n-1}$ (see e.g \cite{schneider-book}, p. 70):
\begin{eqnarray*}
N(K,x):=\{u\in S^{n-1}|(u,y-x)\leq 0 \mbox{ for any } y\in K\}.
\end{eqnarray*}
It is clear that $N(K,x)$ is non-empty if and only if $x$ belongs to
the boundary of $K$. Now define the normal cycle of $K$ by
$$N(K):=\cup_{x\in K} \{(x,u)|\, u\in N(K,x)\}.$$
It is not hard to see that $N(K)$ is a closed subset of $V\times
S^{n-1}$. Moreover it is locally bi-Lipschitz equivalent to
$\RR^{n-1}$, \footnote{This fact was communicated to me by Joe Fu.
Unfortunately I have no reference to it.} and hence integrating of
smooth differential $(n-1)$-forms on $V\times S^{n-1}$ over $N(K)$
defines a continuous linear functional on such forms (more precisely
$N(K)$ can be considered as an integral $(n-1)$-current). A proof of
the following result can be found in \cite{part3}; it is based on
some geometric measure theory and previous work of Fu
\cite{fu-1}-\cite{fu-3}, \cite{fu-94} and other people
\cite{wintgen}, \cite{zahle} on normal cycles (the references can be
found in \cite{part3}).
\begin{proposition}
Let $\ome$ be an infinitely smooth $(n-1)$-form on $V\times
S^{n-1}$. Then the functional
$$K\mapsto \int_{N(K)} \ome$$
is a continuous valuation on $\ck(V)$. If moreover $\ome$ is
invariant with respect to translations in $V$ then the above
expression is a smooth translation invariant valuation in the sense
of Definition \ref{D:smooth-transl}.
\end{proposition}
Let us denote by $\Ome^{n-1}_{tr}(V\times S^{n-1})$ the space of
infinitely smooth $(n-1)$-forms on $V\times S^{n-1}$ which are
invariant with respect to translations on $V$.
\begin{proposition}[\cite{part1}, Theorem
5.2.1]\label{P:normal-transl} The linear map $\CC\oplus
\Ome^{n-1}_{tr}(V\times S^{n-1})\to Val^{sm}(V)$ given by
$$(a,\ome)\mapsto a\cdot vol(K)+\int_{N(K)}\ome$$
is continuous and onto.
\end{proposition}
The proof of this theorem is based on the observation that the map
in the proposition can be rewritten in metric free terms, such that
it will commute with the action of the full liner group $GL(V)$. The
irreducibility theorem implies that the image of this map is dense
in $Val^{sm}(V)$. The fact that the image is closed follows from a
rather general representation theoretical result due to Casselman
and Wallach which says that any morphism between two
$GL(V)$-representations in Fr\'echet spaces satisfying appropriate
technical conditions, has a closed image (see \cite{part1},Theorem
1.1.5, for a precise statement and references).

The kernel of this map was described be Bernig and Br\"ocker
\cite{bernig-brocker} by a system of differential and integral
equations. Bernig has applied very successfully this description in
classification problems of translation invariant valuations
invariant under various groups transitive in spheres
\cite{bernig-sun}, \cite{bernig-spin7-g2},
\cite{bernig-quaternions}.

\subsection{Valuations invariant under a group.}
\label{Ss:invariant-val} Let $G$ be a compact subgroup of the group
of orthogonal transformations of a Euclidean $n$-dimensional space
$V$. We denote by $Val^G$ the subspace of $Val(V)$ of $G$-invariant
valuations. When $G=SO(n)$ the space $Val^G$ was described by
Hadwiger (see Section \ref{Ss:hadwiger-thm}). There are examples of
other groups, e.g. the unitary group $U(n/2)$,  of particular
interest to integral geometry. In fact whenever the space $Val^G$
turns out to be finite dimensional we may hope to classify it
explicitly in geometric terms and then apply this classification to
integral geometry, for example to obtain generalizations of Crofton
and principal kinematic formulas. The first general result in this
direction is as follows.
\begin{proposition}[\cite{alesker-adv-00}]\label{P:abstract-hadwiger}
Let $G$ be a compact subgroup of the orthogonal group. The space
$Val^G$ is finite dimensional if and only if $G$ acts transitively
on the unit sphere.
\end{proposition}
Recall also that by Proposition \ref{P:smooth-invariant-val} if $G$
acts transitively on the sphere then $Val^G\subset Val^{sm}$. This
equips $Val^G$ with the product. Evidently we have also McMullen's
decomposition
$$Val^G=\oplus_{i=0}^n Val^G_i.$$
Thus $Val^G$ becomes a finite dimensional commutative associative
graded algebra with unit. It satisfies the Poincar\'e duality and
two versions of the hard Lefschetz theorem as in Section
\ref{Ss:hard-lefschetz}. Moreover it was shown by Bernig
\cite{bernig-spin7-g2} that for such $G$ all $G$-invariant
valuations are even. Next $Val^G_1=\CC\cdot V_1$,
$Val^G_{n-1}=\CC\cdot V_{n-1}$ by \cite{alesker-gafa-04}.

In topology there is an explicit classification of compact connected
Lie groups acting transitively and effectively on spheres due to A.
Borel and Montgomery-Samelson. There are 6 infinite lists
$$SO(n), U(n), SU(n), Sp(n), Sp(n)\cdot Sp(1), Sp(n)\cdot U(1),$$
and 3 exceptional groups
$$Spin(7), Spin(9), G_2.$$
Valuations in the case of $SO(n)$ were completely studied by
Hadwiger \cite{hadwiger-book}. The next interesting case is the
unitary group $U(n)$. This case turned out to be more complicated
than $SO(n)$ and in recent years there was a considerable progress
in it. There is a complete geometric classification
\cite{alesker-jdg-03}, the description of the algebra structure
\cite{fu-unitary}, and the principal kinematic formula
\cite{bernig-fu-annals}. This is discussed in detail in Fu's
lectures. For most of the other groups new strong results with
applications to integral geometry were obtained recently by Bernig
in a series of articles. We refer to his survey \cite{bernig-survey}
reporting on the progress.

%\subsubsection{Complex and quaternionic
%analysis.}\label{Sss:quaternionic-analysis}

%\newpage

\section{Valuations on manifolds.}\label{S:valuations-on-manifolds}
The notion of valuation on smooth manifolds was introduced by the
author in \cite{part2}. The goal of this section is to describe this
notion, its properties, and some applications to integral geometry
established in \cite{part1}-\cite{part4}, \cite{part3},
\cite{alesker-integral}, \cite{alesker-bernig}. In particular we
extend the product construction to the setting of valuations on
manifolds and explain its intuitive meaning. This intuitive
interpretation is based on another useful notion of generalized
valuation which establishes an explicit link between the valuations
theory and a better studied notion of constructible functions. The
usefulness of this comparison will be illustrated on several other
examples. Next we introduce operations of pull-back and push-forward
under smooth maps of manifolds in a number of important special
cases generalizing familiar operations on smooth functions,
measures, and constructible functions. All these structures are
eventually used to define a general Radon type transform on
valuations which generalizes the classical Radon transforms on
smooth and constructible functions.

\subsection{Definition of smooth valuations on
manifolds; basic examples.}\label{Ss:definitions-manifolds} The
original approach \cite{part2} to define smooth valuations on smooth
manifolds was rather technical. In these notes we will follow a
different, more direct and actually equivalent approach, which
however might look not very natural and less motivated.

Let $X$ be a smooth manifold of dimension $n$.\footnote{All
manifolds are assumed to be countable at infinity, i.e. presentable
as a union of countably many compact subsets.} We describe a certain
class of finitely additive measures on nice subsets of $X$ (to be
more precise, on compact submanifolds with corners). In our current
approach this class is defined by the explicit construction of
integration of a differential form with respect to the normal cycle.
While in the original approach \cite{part1} this description was a
theorem rather than a definition, it seems to be faster not to
repeat all the intermediate steps leading to it. A reader may take
Proposition \ref{P:normal-transl} above as a possible justification
for the current approach.

\hfill

A submanifold with corners of $X$ is a closed subset $P\subset X$
which is locally diffeomorphic to $\RR_{\geq 0}^i\times \RR^j$ where
$i,j$ are integers (then necessarily $0\leq i+j\leq n$). We denote
by $\cp(X)$ the family of all compact submanifolds with corners.
Basic examples from $\cp(X)$ are compact smooth submanifolds,
possibly with boundary. When $X=\RR^n$ simplices, or more generally
simple polytopes, of any dimension belong to $\cp(X)$; however
non-simplicial polytopes (such as the octahedron in $\RR^3$) do not.

We are going to define the normal cycle of $P\in\cp(X)$. Let $T^*X$
denote the cotangent bundle of $X$. Let $\PP_X$ denote the oriented
projectivization of $T^*X$, namely $\PP_X:=(T^*X\backslash
\underline{0})/\RR_{>0}$ where $\underline{0}$ is the zero-section
of $T^*X$, and $\RR_{>0}$ is the multiplicative group of positive
real numbers acting on $T^*X$ by multiplication on the cotangent
vectors. We call $\PP_X$ the {\itshape cosphere bundle} since if one
fixes a Riemannian metric on $X$ then it induces identification of
$\PP_X$ with the unit (co)tangent bundle.

Let $P\in\cp(X)$. Let $x\in P$ be a point. The {\itshape tangent
cone} to $P$ at $x$ is the subset of the tangent space $T_xX$ of all
$\xi$ such that there exists a $C^1$-smooth curve $\gamma\colon
[0,1]\to P$ such that $\gamma(0)=x,\, \gamma'(0)=\xi$. It it not
hard to see that $T_xP\subset T_xX$ is a closed convex polyhedral
cone. Let $(T_xP)^o$ denote the dual cone, namely
$$(T_xP)^o:=\{\eta\in T_x^*X|\, <\eta,\xi>\leq 0 \mbox{ for any }
\xi\in T_xP\}.$$ This is a closed convex cone in $T^*X$. Now define
the {\itshape normal cycle} of $P$ by
$$N(P):=\cup_{x\in P}\left(((T_xP)^o\backslash\{0\})/\RR_{>0}\right).$$
It is well known (and easy to see) that $N(P)$ is a compact
$n-1$-dimensional submanifold of $\PP_X$ with singularities. Also it
is Legendrian with respect to the canonical contact structure on
$\PP_X$ (though this fact will not be used explicitly in these
lectures).
\begin{remark}
If $X=\RR^n$ and $P\in \cp(\RR^n)$ is convex then this definition of
the normal cycle coincides with the definition of the normal cycle
from Section \ref{Sss:normal-cycle}. Actually the normal cycle can
be defined for other classes of sets: sets of positive reach (which
includes convex compact sets in the case $X=\RR^n$), and subanalytic
sets when $X$ is a real analytic manifold (see Fu \cite{fu-94} which
is based on \cite{fu-1}, \cite{fu-2}, \cite{fu-3}, \cite{fu-90}, and
develops further \cite{federer-curv},\cite{wintgen}, \cite{zahle}).
An essentially equivalent notion of characteristic cycle was
developed in \cite{kashiwara-schapira} for subanalytic sets using a
different approach.
\end{remark}

Below in this article we will assume for simplicity of exposition
that $X$ is oriented; this assumption can be easily removed. The
orientation of $X$ induces an orientation of the normal cycle of
every subset.

\begin{definition}\label{D:smooth-val-mflds}
A map $\phi\colon \cp(X)\to \CC$ is called a smooth valuation if
there exist a measure $\mu$ on $X$ and an $n-1$-form $\ome$ on
$\PP_X$, both infinitely smooth, such that
$$\phi(P)=\mu(P)+\int_{N(P)}\ome$$
for any subset $P\in \cp(X)$.
\end{definition}
\begin{remark}
This definition should be compared with Proposition
\ref{P:normal-transl}. It can be shown that any translation
invariant {\itshape convex} valuation on $\RR^n$ which is smooth in
the sense of Definition \ref{D:smooth-transl} can be naturally
extended to a broader class of sets: to compact sets of positive
reach and also to relatively compact subanalytic subsets of $\RR^n$.
This is done as follows: given a convex valuation $\phi\in
Val^{sm}(\RR^n)$, let us represent it (non-uniquely) in the form
$$\phi(K)=a\cdot vol(K)+\int_{N(K)}\ome,$$
where $\ome$ is a smooth translation invariant form. Then $\phi$ can
be extended by the same formula to any subsets from the above
broader class; this extension is independent of the choice of the
form $\ome$ and the constant $a$.
\end{remark}
It can be shown that every smooth valuation is a finitely additive
functional in some precise sense \cite{part2}.

Let us denote by $V^\infty(X)$ the space of all smooth valuations.
The space $V^\infty(X)$ is the main object of study in what follows.

\begin{example}
(1) Any smooth measure on $X$ is a smooth valuation. Indeed let us
take $\ome=0$ in Definition \ref{D:smooth-val-mflds}.

(2) The Euler characteristic $\chi$ is also a smooth valuation. This
fact is less obvious. In the current approach, it is a reformulation
of a version of the Gauss-Bonnet formula due to Chern
\cite{chern-45} who has constructed $\mu$ and $\ome$ to represent
the Euler characteristic; his construction depends on the choice of
a Riemannian metric on $X$.

(3) The next example is very typical for integral geometry. Let
$X=\CC\PP^n$ be the complex projective space. Let $\grc$ denote the
Grassmannian of all complex projective subspaces of $\CC\PP^n$ of a
fixed complex dimension $k$. It is well known that $\grc$ has a
unique probability measure $dE$ invariant under the group $U(n+1)$.
Consider the functional
$$\phi(P)=\int_{E\in \grc}\chi(P\cap E)\,dE.$$
Then $\phi\in V^\infty(\CC\PP^n)$ (this follows e.g. from Fu
\cite{fu-90}).
\end{example}

\hfill

The space $V^\infty(X)$ is naturally a Fr\'echet space. Indeed it is
a quotient space of the direct sum of Fr\'echet spaces
$\cm^\infty(X)\oplus \Ome^{n-1}(\PP_X)$ by a closed subspace, where
$\cm^\infty(X)$ denotes the space of infinitely smooth measures. The
subspace of pairs $(\mu,\ome)$ representing the zero valuation was
described by Bernig-Br\"ocker \cite{bernig-brocker} in terms of a
system of differential and integral equations.

One can show \cite{part2} that smooth valuations form a sheaf. That
means that:

(1) we have the natural restriction map $V^\infty(U)\to V^\infty(V)$
for any open subsets $V\subset U\subset X$;

(2) given an open covering $\{U_\alp\}$ of an open subset $U$, and
$\phi\in V^\infty(U)$ such that the restriction $\phi|_{U_\alp}$ of
$\phi$ to all $U_\alp$ vanishes, then $\phi=0$;

(3) given an open covering $\{U_\alp\}$ of an open $U$ and
$\phi_\alp\in V^\infty(U_\alp)$ for any $\alp$ such that
$\phi_\alp|_{U\alp\cap U_\beta}=\phi_\beta|_{U\alp\cap U_\beta}$ for
all $\alp,\beta$, then there exists (unique by (2)) $\phi\in
V^\infty(U)$ such that $\phi|_{U_\alp}=\phi_\alp$.

\subsection{Canonical filtration on smooth
valuations.}\label{Ss:filtration} \def\cvx{\cv^\infty_X} The space
of smooth valuations carries a canonical filtration by closed
subspaces. In this section we summarize its main properties without
giving a precise definition for which we refer to \cite{part2}. The
important property of this filtration is that it partly allows to
reduce the study of valuations on manifolds to the more familiar
case of translation invariant convex valuations.

Let us denote by $Val(TX)$ the (infinite dimensional) vector bundle
over $X$ such that its fiber over a point $x\in X$ is equal to the
space $Val^\infty(T_xX)$ of smooth translation invariant convex
valuations on $T_xX$. By McMullen's theorem it has a grading by the
degree of homogeneity:
$Val^\infty(TX)=\oplus_{i=0}^nVal^\infty_i(TX).$
\begin{theorem}\label{T:filtration}
There exists a canonical filtration of $V^\infty(X)$ by closed
subspaces $$V^\infty(X)=W_0\supset W_1\supset\dots\supset W_n, \, \,
n=\dim X,$$ such that the associated graded space $\oplus _{i=0}^n
W_i/W_{i+1}$ is canonically isomorphic to the space of smooth
sections $C^\infty(X,Val^\infty_i(TX))$.
\end{theorem}
\begin{remark}\label{R:filtration-spec}
(1) For $i=n$ the above isomorphism means that $W_n$ coincides with
the space of smooth measures on $X$.

(2) For $i=0$ the above isomorphism means that $W_0/W_1$ is
canonically isomorphic to the space of smooth functions
$C^\infty(X)$. The epimorphism $V^\infty(X)\to C^\infty(X)$ with the
kernel $W_1$ is just the evaluation-on-points map
$$\phi\mapsto [x\mapsto \phi(\{x\})].$$
Thus $W_1$ consists precisely of valuations vanishing on all points.

(3) Actually $U\mapsto W_i(U)$ defines a subsheaf $\cw_i$ of the
sheaf of valuations.
\end{remark}

\subsection{Integration functional.}\label{Ss:integration}
Let $V_c^\infty(X)$ denote the subspace of $V^\infty(X)$ of
compactly supported valuations. (The definition is obvious: a
valuations $\phi$ is said to have a compact support if there exists
a compact subset $A\subset X$ such that the restriction
$\phi|_{X\backslash A}=0$.) Clearly if $X$ is compact then
$V_c^\infty(X)=V^\infty(X)$. $V_c^\infty(X)$ carries a natural
locally convex topology such that the natural imbedding to
$V^\infty(X)$ is continuous (however in general this is not a
Fr\'echet space, but rather a strict inductive limit of Fr\'echet
spaces, see Section 5.1 of \cite{part4}).

The integration functional
$$\int_X\colon V_c^\infty(X)\to \CC$$
is defined by $\int_X\phi:=\phi(X)$.

Formally speaking, $\phi(X)$ is not defined when $X$ is not compact.
The formal way to define it is to choose first a large compact
domain $A$ containing the support of $\phi$ and set
$\int_X\phi:=\phi(A)$. Then one can show that this definition is
independent of a large subset $A$. Moreover $\int_X$ is a continuous
linear functional.

\subsection{Product on smooth valuations on manifolds and Poincar\'e
duality.}\label{Ss:product-manifolds} The product on smooth
translation invariant convex valuations which was discussed in
Section \ref{Ss:product}, can be extended to the case of smooth
valuations on manifolds. We will describe below its main properties,
and in Section \ref{Ss:product-generalized} we will explain its
intuitive meaning. However we present no construction of it in these
notes. For the moment there are two different constructions of the
product, both are rather technical. The first one was done in
several steps. Initially the product was constructed by the author
\cite{part1} on $\RR^n$ (earlier the same construction was done even
in a more specific situation \cite{alesker-gafa-04} of convex
valuations polynomial with respect to translations). Then this
construction was extended by Fu and the author \cite{part3} to any
smooth manifold: it was shown that the product can be defined
locally by a choice of a diffeomorphism of $X$ with $\RR^n$ and
applying the above construction, and the main technical point was to
show that the product is independent of the choice of this local
diffeomorphisms. The second and rather different construction of the
product was given recently by Bernig and the author
\cite{alesker-bernig}. This construction describes the product of
valuations directly in terms of the forms $\mu,\ome$ defining the
valuations; it uses the Rumin operator and some other standard
operations on differential forms. In comparison to the first
construction, the second one has the advantage of being independent
of extra structures on $X$ (such as a coordinate system) and also
some other technical advantages. However it is less intuitive than
the first one. We summarize basic properties of the product as the
following.

\begin{theorem}
There exists a canonical product $V^\infty(X)\times V^\infty(X)\to
V^\infty(X)$ which is

(1) continuous;

(2) commutative and associative;

(3) the filtration $W_\bullet$ is compatible with it:
$$W_i\cdot W_j\subset W_{i+j}$$
where we set $W_k=0$ for $k>n=\dim X$;

(4) the Euler characteristic $\chi$ is the unit in the algebra
$V^\infty(X)$;

(5) this product commutes with restrictions to open and closed
submanifolds.

Thus $V^\infty$ is a commutative associative filtered unital algebra
over $\CC$.
\end{theorem}
Let us also add that the evaluation-on-points map $V^\infty(X)\to
C^\infty(X)$ defined in Remark \ref{R:filtration-spec}(2) is an
epimorphism of algebras when $C^\infty(X)$ is equipped with the
usual pointwise product.

\hfill

An important property of the product is a version of the Poincar\'e
duality. Consider the bilinear map
$$V^\infty(X)\times V^\infty_c(X)\to \CC$$
defined by $(\phi, \psi)\mapsto \int_X\phi\cdot \psi$.
\begin{theorem}\label{T:selfduality}
This bilinear form is a perfect pairing. In other words, the induced
map
$$V^\infty(X)\to (V^\infty_c(X))^*$$
is injective and has a dense image with respect to the weak
topology.
\end{theorem}

\subsection{Generalized valuations and constructible functions.}\label{Ss:generalized}
\begin{definition}
The space of generalized valuations is defined by
$$V^{-\infty}(X):=(V^\infty_c(X))^*$$
equipped with the weak topology. Elements of this space are called
generalized valuations.
\end{definition}
By Theorem \ref{T:selfduality} we have the canonical imbedding with
dense image
$$V^\infty(X)\inj V^{-\infty}(X).$$
Informally speaking, at least when $X$ is compact, the space of
valuations is essentially self-dual (up to completion). This
imbedding also means that $V^{-\infty}(X)$ is a completion of
$V^\infty(X)$ in the weak topology. Every smooth valuation can be
considered as a generalized one.

The advantage of working with generalized valuations is that they
contain constructible functions (described below) as a dense
subspace. This gives a completely different point of view on
valuations which is often useful especially on a heuristic level.
Constructible functions have been studies quite extensively by the
methods of algebraic topology (sheaf theory, see the book
\cite{kashiwara-schapira}, Ch. 9). We will illustrate this below
while discussing again the product on valuations, a Radon type
transform, and the Euler-Verdier involution.

\hfill

Let us define the space of constructible functions on $X$. In the
literature there are various slightly different definitions of this
notion, but the differences are technical rather than conceptual.
For simplicity of the exposition we will assume in these notes,
while talking about constructible functions, that $X$ is a real
analytic manifold.
\begin{definition}\label{D:constructible}
A function $f\colon X\to \CC$ on a real analytic manifold $X$ is
called constructible if it takes finitely many values, and for any
$a\in \CC$ the level set $f^{-1}(a)$ is subanalytic.
\end{definition}
For the definition of a subanalytic set see Section 1.2 of
\cite{part4}, or for more details \S 8.2 of the book
\cite{kashiwara-schapira}. Constructible functions with compact
support form a linear space which will be denoted by $\cf(X)$.
Moreover it is an algebra with pointwise product.

An important property of constructible functions is that they also
admit a normal cycle such that if $P\in \cp(X)$ is subanalytic then
the normal cycle of the indicator function $\One_P$ is equal to the
normal cycle of $P$ (see \cite{fu-94}, Ch. 9 of
\cite{kashiwara-schapira}). Using this notion we define the map
$$\Xi\colon \cf(X)\to V^{-\infty}(X)$$
as follows. Let $\phi\in V^{\infty}_c(X)$ be given by
$\phi(P)=\mu(P)+\int_{N(P)}\ome$ with smooth $\mu,\ome$. Then for
any $f\in \cf(X)$
$$<\Xi(f),\phi>=\int_X f\cdot d\mu+ \int_{N(f)}\ome.$$
The map $\Xi$ is well defined, i.e. it is independent of a
particular choice of $\mu,\ome$ representing $\phi$. $\Xi$ is a
linear injective map with dense image (\cite{part4}, Section 8.1).

\hfill

To summarize, we have a large space of generalized valuations with
two completely different dense subspaces
\begin{eqnarray}\label{E:subspaces-valuations}
V^\infty(X)\subset V^{-\infty}(X)\supset \cf(X).
\end{eqnarray}
Notice that when $X$ is compact, the image of the constant function
$1\in \cf(X)$ in $V^{-\infty}(X)$ coincides with the image of the
Euler characteristic $\chi\in V^\infty(X)$.

While working with valuations it is useful to keep in mind the
imbeddings (\ref{E:subspaces-valuations}). The role of constructible
functions in the theory of valuations is somewhat analogous to the
role of delta-functions in the classical theory of generalized
functions (distributions). Often it is instructive to compare
various structures on valuations with their analogues on
constructible functions. We will see several examples of this below.
Here we will show how it works for the integration functional and
the filtration $W_\bullet$.

It was shown in \cite{part4} that the integration functional
$\int_X\colon V^\infty_c\to \CC$ extends uniquely by continuity in
weak topology to the generalized valuations with compact support:
$$\int_X\colon V^{-\infty}_c(X)\to \CC.$$
Let us restrict this functional to the subspace $\cf_c(X)$ of
constructible functions with compact support. It turns out that this
restriction coincides with the integration with respect to the Euler
characteristic; this operation is uniquely characterized by the
property that for any compact subanalytic subset $P\subset X$
$$\int_X\One_P=\chi(P).$$

Let us consider now the filtration $W_\bullet$ on $V^\infty(X)$. Let
$W_i'$ denote the closure of $W_i$ in $V^{-\infty}(X)$ with respect
to the weak topology. By \cite{part4} the restriction of $W_i'$ back
to $V^\infty(X)$ coincides with $W_i$: $W_i'\cap V^\infty(X)=W_i$.
Consider the induced filtration on constructible functions, namely
$$\cf(X)=\cf(X)\cap W_0'\supset \cf(X)\cap W_1'\supset\dots\supset \cf(X)\cap
W_n'.$$ It was shown in \cite{part4} that $\cf(X)\cap W_i'$ consists
of constructible functions with codimension of the support at most
$i$. In particular $\cf(X)\cap W_n'$ consists of functions with
discrete support.

\subsection{Euler-Verdier involution.}\label{Ss:euler-verdier}
Let us give another example of an application of the comparison with
constructible functions. The space of constructible functions has a
canonical linear involution called the Verdier involution (see e.g.
\cite{kashiwara-schapira}). In the special case of functions on
$\RR^n$ which are constructible in a more narrow (polyhedral) sense
this involution has been known for convexity experts under the name
of Euler involution. We will see that it extends naturally to
valuations, and this extension will be called the Euler-Verdier
involution.

Here we will choose a sign normalization different from the standard
one. Let us describe the Verdier involution $\sigma$ (with a
different sign convention) in a special case when a constructible
function has the form $\One_P$ where $P$ is a compact subanalytic
submanifold with corners (the general case is not very far from this
one using the linearity property of it). Then
$$\sigma(\One_P)=(-1)^{n-\dim P}\One_{\textrm{int}\, P},$$
where $\textrm{int}\, P$ denotes the relative interior of $P$. One
has $\sigma^2=Id$.

\begin{theorem}[\cite{part4}]\label{T:euler-verdier}

(1) The involution $\sigma$ extends (uniquely) by continuity to
$V^{-\infty}(X)$ in the weak topology. This extension is also
denoted by $\sigma$;

(2) $\sigma$ preserves the class of smooth valuations and
$\sigma\colon V^\infty(X)\to V^\infty(X)$ is a continuous linear
operator (in the Fr\'echet topology);

(3) $\sigma^2=Id$;

(4) $\sigma\colon V^\infty(X)\to V^\infty(X)$ is an algebra
automorphism;

(5) $\sigma$ preserves the filtration $W_\bullet$, namely
$\sigma(W_i)=W_i$;

(6) for any smooth translation invariant valuation $\phi$ on $\RR^n$
one has $$(\sigma \phi)(K)=(-1)^{\deg \phi}\phi(-K)$$ where $\deg
\phi$ denotes the degree of homogeneity of $\phi$;

(7) $\sigma$ commutes with restrictions to open subsets (both for
smooth and generalized valuations).
\end{theorem}
\begin{remark}
Though $\sigma$ was defined above only on a real analytic manifold
$X$, it can be defined on any smooth manifold as a continuous linear
operator $\sigma\colon V^{-\infty}(X)\to V^{-\infty}(X)$. Then it
satisfies the properties (2)-(7) of Theorem \ref{T:euler-verdier}.
\end{remark}

\subsection{Partial product on generalized
valuations.}\label{Ss:product-generalized} In this section we
discuss the promised intuitive meaning of the product on valuations.
This interpretation was conjectured by the author
\cite{alesker-survey} and proved rigorously by Bernig and the author
\cite{alesker-bernig}. It provides yet another example of the
relevance of constructible functions to valuations.

Recall again that we have the imbedding of smooth valuations to
generalized ones:
$$V^\infty(X)\subset V^{-\infty}(X).$$
One could try to extend the product on smooth valuations to
$V^{-\infty}(X)$ say by continuity. Unfortunately this is not
possible. The situation here is much analogous to what is known in
the classical theory of generalized functions (see e.g.
\cite{hormander-book}). There the space of smooth functions
$C^\infty(X)$ is naturally imbedded to the larger space of
generalized functions $C^{-\infty}(X)$ which is a completion of it
in the weak topology. The space $C^\infty(X)$ has its usual
pointwise product. However this product does not extend to
$C^{-\infty}(X)$ by continuity: for example no rigorous way is known
how to take the square of the delta-function on $X=\RR$.
Nevertheless it is still possible to define a {\itshape partial}
product on $C^{-\infty}(X)$. That roughly means that one can define
a product of two generalized functions whose "singularities" are
disjoint. The precise technical condition is formulated in the
language of the wave front sets of generalized functions in the
sense of H\"ormander-Sato; we will not reproduce it here, but rather
refer to \cite{hormander-book}. This partial product is natural and
satisfies some continuity properties (\cite{guillemin-sternberg},
Ch. VI \S 3).

In the case of valuations we have the following result.
\begin{theorem}[\cite{alesker-bernig}]
There exists a partial product on $V^{-\infty}(X)$ extending the
product on $V^\infty(X)$. It is commutative and associative.
\end{theorem}
We refer to \cite{alesker-bernig} for the precise technical
formulation when the partial product of two generalized valuations
is defined. We notice only that the condition is also formulated in
the language of wave front sets.

\hfill

Now we can try to restrict the partial product on generalized
valuations to constructible functions and see what we get. The
answer is very natural: we just get their pointwise product (under
certain technical assumptions on the functions guaranteeing that
their product in $V^{-\infty}(X)$ is well defined). More precisely
we have the following result.
\begin{theorem}[\cite{alesker-bernig}]\label{T:product-constr}
Let $P,Q\subset X$ be compact submanifolds with corners which
intersect transversally. Then the product of $\One_P$ and $\One_Q$
in the sense of generalized valuations is well defined and is equal
to $\One_{P\cap Q}$ (notice that under the transversality
assumption, $P\cap Q$ is also a compact submanifold with corners).
\end{theorem}

We did not give a formal definition of transversality of two
submanifolds with corners. In the special case of submanifolds
without corners, the definition is the usual one. In general the
precise definition is given in \cite{alesker-bernig}. Notice only
that any two compact submanifolds with corners can be made
transversal to each other by applying to one of them a generic
diffeomorphism which is arbitrarily close (in the
$C^\infty$-topology) to the identity map.

\subsection{Few examples of computation of the product in integral
geometry.}\label{Ss:product-IG} In this section we give few examples
of computation of the product of valuations in the complex
projective space $\CC\PP^n$. These examples are very typical in
integral geometry. We present a heuristic argument since hopefully
it will clarify the intuition behind the product in applications.

First let us give few heuristic remarks of a general nature. Let $P$
be a compact real analytic subset of a real analytic manifold $X$.
As we have seen in Section \ref{Ss:generalized}, the indicator
function $\One_P$ can be considered as a generalized valuation. Can
we consider this non-smooth valuation as a finitely additive
measure? The answer is: "essentially yes". This measure is defined
only on nice compact subsets of $X$ which are "transversal" to $P$.
It is equal to
$$K\mapsto \chi(K\cap P).$$
We do not want to formalize this here, but this is the right way to
think of $\One_P$ as a measure.

\hfill

\def\pgc{{}\!^ {\textbf{C}} G}

Let now $X=\CC\PP^n$ with the Fubini-Study metric. Let us denote by
$\pgc_l$ the Grassmannian of $l$-dimensional complex projective
subspaces of $\CC\PP^n$. Clearly it is equal to the Grassmannian of
$l+1$-dimensional complex linear subspaces in $\CC^{n+1}$. Let us
consider the smooth $U(n+1)$-invariant valuations
\begin{eqnarray}\label{E:phi-l}
\phi_l(K):=\int_{\pgc_l} \chi(K\cap E) dE
\end{eqnarray}
where $dE$ is the Haar measure on $\pgc$  normalized somehow (we do
not care about normalization constants). We claim that
\begin{eqnarray}\label{E:product-example1}
\phi_l\cdot \phi_m=\left\{\begin{array}{cc}
                          c\cdot \phi_{l+m-n},& l+m\geq n\\
                          0,& l+m<n
                          \end{array}\right.
\end{eqnarray}
where $c\ne 0$ is a normalizing constant depending on normalizations
of Haar measures.

Let us give a heuristic proof of this equality. Using the discussion
in the beginning of this section, we observe that
\begin{eqnarray*}
\phi_l(K)=\left(\int_{\pgc_l}\One_E dE\right)(K)
\end{eqnarray*}
where $\One_E$ is considered as a generalized valuation. Hence
\begin{eqnarray*}
\phi_l\cdot\phi_m=\int_{(E,F)\in \pgc_l\times\pgc_m}\One_E\cdot
\One_F \,dE\, dF=\int_{\pgc_l\times\pgc_m} \One_{E\cap F}\, dE\,dF
\end{eqnarray*}
where the last equality is due to Theorem \ref{T:product-constr}.
Since for generic projective subspaces $E$ and $F$ their
intersection $E\cap F$ is a projective subspace of dimension $l+m-n$
for $l+m\geq n$ and empty otherwise, it follows that
$$\int_{\pgc_l\times \pgc_m}\One_{E\cap F}\,dE\, dF=c
\int_{\pgc_{l+m-n}}\One_M dM=c\cdot \phi_{l+m-n}.$$ Thus the
equality (\ref{E:product-example1}) is proved.

\hfill

Let us consider another important example of the product on
$\CC\PP^n$. We claim that the $U(n+1)$-invariant valuation
\begin{eqnarray}\label{E:phi-2}
K\mapsto \int_{\pgc_l}V_i(K\cap E) dE
\end{eqnarray}
is equal to $\phi_l\cdot V_i$ where $\phi_l$ is defined by
(\ref{E:phi-l}). First observe
$$V_i(K\cap E)=\One_{K\cap E}(V_i)=\int \One_{K\cap E}\cdot V_i$$
where $\One_{K\cap E}$ is considered as a generalized valuation,
$\int$ in the last expression is the integration functional (i.e the
evaluation on the whole space $\CC\PP^n$); here both equalities are
tautological by unraveling the definitions. Now we again use Theorem
\ref{T:product-constr} to write (under transversality assumptions)
$\One_{K\cap E}= \One_K\cdot \One_E$. Thus
$$\int \One_{K\cap E}\cdot V_i=\int\One_K\cdot \One_E\cdot
V_i=(\One_E\cdot V_i)(K).$$ Thus the valuation (\ref{E:phi-2}) is
equal to
\begin{eqnarray*}
\int_{\pgc_l}\One_E\cdot V_i \,dE=\left(\int_{\pgc_l}\One_E\,
dE\right)\cdot V_i=\phi_l\cdot V_i
\end{eqnarray*}
as it was claimed.

\hfill

Finally let us compute a generalization of the two previous
examples. We claim
\begin{eqnarray}\label{E:general-example}
\left(\int_{\pgc_l}V_i(\bullet \cap E) dE\right)\cdot
\left(\int_{\pgc_m}V_j(\bullet\cap F) dF\right)=c' \cdot
\int_{\pgc_{l+m-n}}V_{i+j}(\bullet\cap M)dM
\end{eqnarray}
where $c'$ is a constant which can be computed explicitly. By the
previous two examples of this section, Example
\ref{Ex:product-intrinsic} from Section \ref{Ss:product}, and using
the associativity and the commutativity of the product, we see that
the left hand side of (\ref{E:general-example}) is equal to
\begin{eqnarray*}
(\phi_l\cdot V_i)\cdot (\phi_m\cdot V_j)=(\phi_l\cdot \phi_m)\cdot
(V_i\cdot V_j)=c' \cdot \phi_{l+m-n}\cdot V_{i+j}=\mbox{r.h.s. of
(\ref{E:general-example})}.
\end{eqnarray*}
Thus (\ref{E:general-example}) is proved. \qed

\subsection{Functorial properties of
valuations.}\label{Ss:functorial} We describe the operations of
pull-back and push-forward on valuations under smooth maps of
manifolds. These operations generalize the well known operation of
pull-back on smooth and constructible functions, the operation of
push-forward on measures, and integration with respect to the Euler
characteristic along the fibers (also called push-forward) on
constructible functions. However for the moment this is done
rigorously only in several special cases of maps (say submersions
and immersions). We believe however than these constructions can be
extended to "generic" smooth maps as partially defined maps on
valuations. The precise conditions under which the maps could be
defined might be rather technical. For this reason we describe first
the general picture heuristically. This picture should be considered
as conjectural. Then we formulate several rigorous results with
precise conditions under which one can define pull-back and
push-forward on valuations. These special cases turn out to be
sufficient to define rigorously the Radon type transform on
valuations (again under some conditions) in the next section. The
results of this section have been obtained by the author in
\cite{alesker-integral}.

\hfill

Let us start with the heuristic picture. Let us denote by $V(X)$ a
space of valuations on a manifold $X$ without specifying exactly the
class of smoothness (smooth, generalized, or something else).
$V_c(X)$ denotes the subspace of $V(X)$ of compactly supported
valuations. Let $f\colon X\to Y$ be a smooth map of manifolds. There
should exist a partially defined linear map, called push-forward,
$$f_*\colon V_c(X)\dashrightarrow V_c(Y)$$
such that for any nice subset $P\subset Y$
\begin{eqnarray}\label{E:push}
(f_*\phi)(P)=\phi(f^{-1}(P)).
\end{eqnarray}
Since (smooth) measures are contained in $V(X)$, we can make their
push-forward in the sense of valuations. Clearly this operation
should coincide with the classical push-forward of measures.

It immediately follows from (\ref{E:push}) that for the composition
of maps we should have
\begin{eqnarray}\label{E:composition}
(f\circ g)_*=f_*\circ g_*.
\end{eqnarray}

We expect that the following interesting property of push-forward
$f_*$ holds. It should extend somehow to a partially defined map on
generalized valuations. Hence then $f_*$ can be restricted to a
partially defined map on constructible functions; it should be
defined on constructible functions which are "in generic position"
to the map $f\colon X\to Y$. We expect that when $f$ is a proper map
(i.e. preimage of any compact set is compact) then on constructible
functions $f_*$ coincides with the integration with respect to the
Euler characteristic along the fibers. Let us remind this operation
assuming that $X$ and $Y$ are real analytic manifolds and $f$ is a
proper real analytic map. It is uniquely characterized by the
following property: Let $P\subset X$ be a subanalytic compact
subset. Then $(f_*\One_P)(y)=\chi(P\cap f^{-1}(y))$ for any point
$y\in Y$. One can show that $f_*$ maps constructible functions to
constructible ones. We refer to Ch. 9 of \cite{kashiwara-schapira}
for further details.

The push-forward should be related to the filtration on valuations
in the following way
$$f_*(W_i)\subset W_{i-\dim X+\dim Y}.$$
Also $f_*$ should commute (up to a sign) with the Euler-Verdier
involution.

\hfill

Let us now discuss the pull-back operation
$$f^*\colon V(Y)\dashrightarrow V(X)$$
which should be a partially defined linear map in the opposite
direction. Heuristically, $f^*$ should be the dual map to $f_*$
(recall from Section \ref{Ss:product-manifolds} that $V_c(X)$ and
$V(X)$ are essentially dual to each other). The pull-back $f^*$
should be a homomorphism of algebras of valuations (again, the
product might be partially defined). We expect $f^*\chi=\chi$. Also
$f^*$ should preserve the filtration
$$f^*(W_i)\subset W_i,$$ and $f^*$ should commute with the
Euler-Verdier involution. Notice that since in particular
$f^*(W_1)\subset W_1$, $f^*$ induces a map between the quotients
$$f^*\colon V(Y)/W_1\to V(X)/W_1.$$
But by Remark \ref{R:filtration-spec}(2), $V(Y)/W_1$ coincides with
functions on $Y$ of an appropriate class of smoothness. In
particular we should get a map
$$f^*\colon C^\infty(Y)\to C^{-\infty}(X).$$
We expect that this is the usual pull-back on smooth functions, i.e.
\begin{eqnarray}\label{E:pull-functions}
f^*(F)=F\circ f.
\end{eqnarray}
Now let us restrict $f^*$ to constructible functions. We expect that
it coincides again with the usual pull-back on constructible
functions which is defined by the same formula
(\ref{E:pull-functions}).

Finally consider the restriction of $f^*$ to (say smooth) measures
on $Y$. In the classical measure theory the operation of pull-back
of a measure does not exist. Nevertheless it is possible to define
it as a valuation, at least under appropriate technical conditions
on the map $f$. Let $\mu$ be a smooth measure on  $Y$. Then, leaving
all the technicalities aside, one should have
$$(f^*\mu)(P)=\int_{y\in Y}\chi(P\cap f^{-1}(y))d\mu(y).$$
In particular if $f\colon X\to Y$ is a linear projection of vector
spaces, and $P\subset X$ is a convex compact subset then
$(f^*\mu)(P)=\mu(f(P))$ is the measure of the projection of $P$.

\hfill

Now let us describe several rigorous results which will be used
later. Let $f\colon X\to Y$ be a smooth map.

\underline{Case 1.} Assume that $f$ is a closed imbedding. Then the
obvious restriction map $V^\infty(Y)\to V^{\infty}(X)$ defines the
pull-back map $f^*$ which is a linear continuous operator. Dualizing
it, we get the push-forward map
$$f_*\colon V^{-\infty}(X)\to V^{-\infty}(Y)$$
which is a linear continuous (in the weak topology) operator. Notice
that in this situation $f_*$ does not preserve the class of smooth
valuations.

It was shown in \cite{alesker-integral} that in this case
$f_*(\One_P)=\One_{f(P)}$ for any compact submanifold with corners
$P\subset X$. It was also shown that $f^*$ can be extended to a
partially defined map $V^{-\infty}(Y)\dashrightarrow V^{-\infty}(X)$
such that if $Q\subset Y$ is a compact submanifold with corners
which is transversal to $X$ then $f^*\One_Q$ is well defined in the
sense of valuations and is equal to $\One_{X\cap Q}$, i.e. the
pull-back on valuations is compatible with the pull-back on
constructible functions.

\underline{Case 2.} Assume that $f$ is a proper submersion. Let us
define the pull-back $f_*\colon V^\infty_c(X)\to V^\infty_c(Y)$ by
$(f_*\phi)(P)=\phi(f^{-1}(P))$ for any compact submanifold with
corners $P\subset Y$. Notice that in this case $f^{-1}(P)$ is a
compact submanifold with corners, and $f_*\phi$ is indeed a smooth
valuation. The map constructed is linear and continuous. Taking the
dual map, we define the pull-back map
$$f^*\colon V^{-\infty}(Y)\to V^{-\infty}(X).$$ It was shown in
\cite{alesker-integral} that in this case for any compact
submanifold with corners $P\subset Y$ one has
$f^*(\One_P)=\One_P\circ f=\One_{f^{-1}(P)}$. It was also shown that
the push-forward $f_*$ extends to a partially defined map on
generalized valuations. However its compatibility with the
integration with respect to the Euler characteristic along the
fibers was proved only under rather ugly restrictions on the class
of constructible functions.

\subsection{Radon transform on valuations on
manifolds.}\label{Ss:radon} In this section we combine the product,
pull-back, and push-forward on valuations to define a Radon type
transform on them. Before we introduce this notion, it is
instructive to remind the general Radon transform on smooth
functions following Gelfand, and less classical but still known
Radon transform on constructible functions. These two completely
different transforms can be considered as special cases of the
general Radon transform on valuations. In our opinion, this is the
most interesting property of the new Radon transform on valuations.

\begin{definition}
A double fibration is a triple of smooth manifolds $X,Y,Z$ with two
submersive maps
$$X\overset{p}{\leftarrow}Z\overset{q}{\to}Y$$
such that the map $Z\overset{p\times q}{\to}X\times Y$ is a closed
imbedding.
\end{definition}

To define a general Radon transform on smooth functions let us fix a
double fibration as above and an infinitely smooth measure $\gamma$
on $Z$. Let us also assume that $q\colon Z\to Y$ is proper. The
Radon transform is the operator $R_\gamma\colon C^\infty_c(X)\to
\cm^\infty(Y)$ (where $\cm^\infty(Y)$ denotes the space of smooth
measures) defined by
\begin{eqnarray}\label{E:radon-class}
R_\gamma f:=q_*(\gamma\cdot p^*f),
\end{eqnarray}
where $p^*f=f\circ p$ is the usual pull-back on smooth functions,
the product is just the usual product of a measure by a function,
and $q_*$ is the usual push-forward on measures. Notice that all
classical Radon transforms on smooth functions have such a form. For
example let us take $X=\RR^n$, $Y$ the Grassmannian of affine
$k$-dimensional subspace, and $Z$ the incidence variety, i.e.
$Z=\{(x,E)\in X\times Y|\, x\in E\}$. Let $\gamma$ be a Haar measure
on $Z$ invariant under the group of all isometries of $\RR^n$. Then
$R_\gamma$ is the classical Radon transform given by integration of
a function on $\RR^n$ over all affine $k$-dimensional subspaces.
There is a very extensive literature on this subject, see e.g.
\cite{gelfand-graev-rosu}, \cite{gelfand-graev-schapiro},
\cite{helgason}.

\hfill

Let us recall the Radon transform with respect to the Euler
characteristic on constructible functions. It was studied for the
real projective spaces and a somewhat restrictive class of
constructible functions by Khovanskii and Pukhlikov
\cite{khovanskii-pukhlikov}; their work has been motivated by the
earlier work of Viro \cite{viro} on Radon transform on complex
constructible functions on the complex projective spaces. We will
discuss and generalize the Khovanskii-Pukhlikov result in the next
section. For subanalytic constructible functions and other spaces
the Radon transform with respect to the Euler characteristic was
studied by Schapira \cite{schapira}. Thus let
$$X\overset{p}{\leftarrow}Z\overset{q}{\to}Y$$
be a double fibration of real analytic spaces with real analytic
maps $p,q$. We assume again that $q$ is proper. Let us denote by
$\cf(X)$ the space of constructible functions as defined in Section
\ref{Ss:generalized}. One defines the Radon transform $R\colon
\cf(X)\to \cf(Y)$ by
\begin{eqnarray}\label{E:radon-constr}
Rf:=q_*p^*(f), \end{eqnarray} where $p^*$ denotes the usual
pull-back on (constructible) functions, and $q_*$ is the integration
with respect to the Euler characteristic along the fibers of $q$.

\hfill

With these preliminaries let us introduce the Radon transform on
valuations. We fix a double fibration as above with the map $q$
being proper. Let us fix a smooth valuation $\gamma\in V^\infty(Z)$.
Define the Radon transform on valuations $\car_\gamma\colon
V^\infty(X)\to V^{-\infty}(Y)$ by
$$\car_\gamma(\phi)=q_*(\gamma\cdot p^*\phi),$$ where $p^*$ and $q_*$ are the
pull-back and push-forward on valuations respectively, and the
product with $\gamma$ is taken in the sense of valuations. It was
shown in \cite{alesker-integral} that the operator $\car_\gamma$ is
a well defined continuous linear operator.

Let us comment on some of the technical difficulties in this
construction. Usually $p^*\phi$ is not a smooth valuation, though
$\phi$ is. Thus we have to multiply the smooth valuation $\gamma$ by
the non-smooth $p^*\phi$. This is always possible in the class of
generalized valuations, but the product is not a smooth valuation.
Next we have to take the push-forward of this generalized valuation.
The push-forward of a generalized valuation under a general proper
submersion is not always defined, but only under some rather
technical condition of "generic position" of "singularities" of the
valuation with respect to the map $q$. Fortunately this technical
condition is satisfied for valuations of the form $\gamma\cdot
p^*\phi$ with smooth $\phi$. It was also shown in
\cite{alesker-integral} that under extra assumptions on the double
fibration the image $\car_\gamma(V^\infty(X))$ is contained in
smooth valuations. Also under a similar extra assumption
$\car_\gamma$ can be extended uniquely by continuity in the weak
topology to generalized valuations $V^{-\infty}(X)$. An example
satisfying both assuptions will be considered in the next section.

\hfill

Let us discuss now the relation of the new Radon transform on
valuations to the classical Radon transforms discussed above in this
section. First let us assume that the valuation $\gamma\in
V^\infty(Z)$ is in fact a smooth measure considered as a smooth
valuation. Then the Radon transform
$$\car_\gamma\colon V^\infty(X)\to V^{-\infty}(Y)$$
vanishes on $W_1\subset V^\infty(X)$. Indeed $p^*(W_1)\subset W_1$,
and $\gamma\cdot W_1=0$ since $\gamma$ is a measure. Hence
$\car_\gamma$ factorizes (uniquely) via the quotient
$V^\infty(X)/W_1=C^\infty(X)$. Notice also that in this case
$\car_\gamma$ takes values in measures, in fact in infinitely smooth
ones. Hence we get a map $C^\infty(X)\to \cm^\infty(Y)$. It was
shown in \cite{alesker-integral} that this map coincides with the
classical Radon transform $R_\gamma$ defined by
(\ref{E:radon-class}).

Let us consider another extremal case of $\car_\gamma$ with
$\gamma=\chi$ being the Euler characteristic. In this case our
discussion will be less rigorous. First assume that $\car_\gamma$
extends naturally to a partially defined map on generalized
valuations $V^{-\infty}(X)\dashrightarrow V^{-\infty}(Y)$. We expect
that its restriction to the class of constructible functions
coincides with the Radon transform with respect to the  Euler
characteristic defined previously by (\ref{E:radon-constr}). This
result was proved rigorously in \cite{alesker-integral} in very
special circumstances. It is desirable to make the result rigorous
under more general assumptions.

\subsection{Khovanskii-Pukhlikov type inversion formula for the Radon
transform on valuations on $\RR\PP^n$.}\label{Ss:inversion} Let us
consider the Radon type transform on valuations in the following
special case. Let $X=\RR\PP^n$ be the real projective space, i.e.
the manifold of lines in $\RR^{n+1}$ passing through the origin. Let
$Y=\RR\PP^{n\vee}$ be the dual projective space, i.e. the manifold
of linear hyperplanes in $\RR^{n+1}$. Let $Z\subset X\times Y$ be
the incidence variety
$$Z:=\{(l,E)\in \RR\PP^n\times\RR\PP^{n\vee}|\, l\subset E\}.$$ We
have the double fibration
$$\RR\PP^n\overset{p}{\leftarrow}Z\overset{q}{\to}\RR\PP^{n\vee}$$
where $p,q$ are the obvious projections. All the manifolds and maps
are real analytic.

We consider the Radon transform
$$\car_\chi\colon V^\infty(\RR\PP^n)\to V^{-\infty}(\RR\PP^{n\vee})$$
with the kernel $\gamma=\chi$ being the Euler characteristic on $Z$.
In this case
$$\car_\chi=q_*p^*.$$ It was shown in
\cite{alesker-integral} that the image of this transformation is
contained in smooth valuations, and $\car_\chi\colon
V^\infty(\RR\PP^n)\to V^\infty(\RR\PP^{n\vee})$ is continuous.
Moreover this operator extends (uniquely) to a continuous linear
operator, also denoted by $\car_\chi$, on generalized valuations
equipped, as usual, with the weak topology:
$$\car_\chi\colon V^{-\infty}(\RR\PP^n)\to
V^{-\infty}(\RR\PP^{n\vee}).$$ It was shown in
\cite{alesker-integral} that $\car_\chi$ is invertible for odd $n$,
and for even $n$ its kernel consists precisely of multiples of the
Euler characteristic. In both cases there is an explicit inversion
formula (in the latter case, up to a multiple of the Euler
characteristic); it generalizes and was motivated by the
Khovanskii-Pukhlikov inversion formula for constructible functions
\cite{khovanskii-pukhlikov}. In order to state the result let us
consider the analogous operator in the opposite direction
$$\car^t_\chi\colon V^{-\infty}(\RR\PP^{n\vee})\to
V^{-\infty}(\RR\PP^n),$$ namely $$\car^t_\chi:=p_*q^*.$$
\begin{theorem}[\cite{alesker-integral}]\label{T:radon-val}
For any generalized valuation $\phi\in V^{-\infty}(\RR\PP^n)$ one
has
\begin{eqnarray}\label{inversion-radon-val-constr}
(-1)^{n-1}\car_\chi^t\car_\chi(\phi)=\phi+\frac{1}{2}((-1)^{n-1}-1)
\left(\int_{\RR\PP^n}\phi\right)\cdot \chi.
\end{eqnarray}
%In particular if $n$ is odd then the Radon transform $\car_\chi$ is
%injective, and (\ref{inversion-radon-val-constr}) is the inversion
%formula. If $n$ is even then the kernel of $\car_\chi$ consists
%precisely of the multiples of the Euler characteristic $\chi$.
\end{theorem}
Let us say a few words on the proof of this theorem. After all the
operators involved were defined, the next technically non-trivial
step was to show that the restriction of $\car_\chi$ to a rather
special class of constructible functions, which is still dense in
$V^{-\infty}(\RR\PP^n)$, coincides with the Radon transform with
respect to the Euler characteristic on constructible functions; also
an analogous result holds for $\car^t_\chi$. Then Theorem
\ref{T:radon-val} follows immediately by continuity from the
Khovanskii-Pukhlikov inversion formula for constructible functions
which claims precisely the identity
(\ref{inversion-radon-val-constr}) for such functions in place of
$\phi$.


\begin{thebibliography}{99}
%\bibitem{alesker-annals}
%Alesker, Semyon; Continuous rotation invariant valuations on convex
%sets. Ann. of Math. (2) 149 (1999), no. 3, 977--1005.
%\bibitem{alesker-geom-dedicata} Alesker, Semyon; Description of
%continuous isometry covariant valuations on convex sets. Geom.
%Dedicata 74 (1999), no. 3, 241--248.

\bibitem{alesker-adv-00}
Alesker, Semyon; On P. McMullen's conjecture on translation
invariant valuations. Adv. Math. 155 (2000), no. 2, 239--263.
\bibitem{alesker-gafa-01}
Alesker, Semyon; Description of translation invariant valuations on
convex sets with solution of P. McMullen's conjecture. Geom. Funct.
Anal. 11 (2001), no. 2, 244--272.
\bibitem{alesker-jdg-03}
Alesker, Semyon; Hard Lefschetz theorem for valuations, complex
integral geometry, and unitarily invariant valuations. J.
Differential Geom. 63 (2003), no. 1, 63--95. math.MG/0209263.
\bibitem{alesker-gafa-04}
Alesker, Semyon; The multiplicative structure on polynomial
continuous valuations. Geom. Funct. Anal. 14 (2004), no. 1, 1--26,
math.MG/0301148.
\bibitem{alesker-su2}
Alesker, Semyon; $SU(2)$-invariant valuations. {\itshape Geometric
aspects of functional analysis}, 21--29, Lecture Notes in Math.,
1850, Springer, Berlin, 2004.
\bibitem{alesker-gafa-sem-04}
Alesker, Semyon; Hard Lefschetz theorem for valuations and related
questions of integral geometry. {\itshape Geometric aspects of
functional analysis}, 9--20, Lecture Notes in Math., 1850, Springer,
Berlin, 2004.
\bibitem{alesker-adv-05}
Alesker, Semyon; Valuations on convex sets, non-commutative
determinants, and pluripotential theory. Adv. Math. 195 (2005), no.
2, 561--595.
\bibitem{alesker-survey-alg-analys}
Alesker, Semyon; Quaternionic plurisubharmonic functions and their
applications to convexity. Algebra i Analiz 19 (2007), no. 1, 5--22;
translation in St. Petersburg Math. J. 19 (2008), no. 1, 1--13.
\bibitem{part1} Alesker, Semyon; Theory of valuations on manifolds, I. Linear spaces. Israel
Journal of Mathematics 156 (2006), 311-339. Also: math.MG/0503397.
\bibitem{part2}
Alesker, Semyon; Theory of valuations on manifolds. II. Adv. Math.
207 (2006), no. 1, 420--454. Also: math.MG/0503399.
\bibitem{part4}
Alesker, Semyon; Theory of valuations on manifolds, IV. New
properties of the multiplicative structure. {\itshape Geometric
aspects of functional analysis}, 1--44, Lecture Notes in Math.,
1910, Springer, Berlin, 2007. Also: math.MG/0511171.
\bibitem{alesker-survey}
Alesker, Semyon; Theory of valuations on manifolds: a survey. Geom.
Funct. Anal. 17 (2007), no. 4, 1321--1341. Also: math.MG/0603372.
\bibitem{alesker-spin9}
Alesker, Semyon; Plurisubharmonic functions on the octonionic plane
and Spin(9)-invariant valuations on convex sets. J. Geom. Anal. 18
(2008), no. 3, 651--686.
\bibitem{alesker-fourier}
Alesker, Semyon; A Fourier type transform on translation invariant
valuations on convex sets. Israel J.Math., to appear. Also:
arXiv:math/0702842.
\bibitem{alesker-integral}
Alesker, Semyon; Valuations on manifolds and integral geometry.
Accepted to Geom. Func. Anal. Also:  arXiv:0905.4046.
\bibitem{alesker-bernig}
Alesker, Semyon; Bernig, Andreas; The product on smooth and
generalized valuations. Accepted to Amer. J. Math. Also:
arXiv:0904.1347.

\bibitem{alesker-bernstein}
Alesker, Semyon; Bernstein, Joseph; Range characterization of the
cosine transform on higher Grassmannians. Adv. Math. 184 (2004), no.
2, 367--379.
\bibitem{part3}
Alesker, Semyon; Fu, Joseph H. G.; Theory of valuations on
manifolds, III. Multiplicative structure in the general case. Trans.
Amer. Math. Soc. 360 (2008), no. 4, 1951--1981. Also:
math.MG/0509512.
\bibitem{bernig-sun}
Bernig, Andreas; A Hadwiger-type theorem for the special unitary
group. Geom.Func.Anal. 19 (2009), 356-372.
\bibitem{bernig-spin7-g2}
Bernig, Andreas; Integral Geometry under $G_2$ and $Spin(7)$. Israel
J. Math., to appear.
\bibitem{bernig-quaternions}
Bernig, Andreas; Invariant valuations on quaternionic vector spaces.
Preprint  arXiv:1005.3654.
\bibitem{bernig-survey}
Bernig, Andreas; Algebraic integral geometry. Preprint
arXiv:1004.3145.
\bibitem{bernig-brocker}
Bernig, Andreas; Br\"ocker, Ludwig; Valuations on manifolds and
Rumin cohomology. J.Differential Geom. 75 (2007), 433-457.
\bibitem{bernig-fu-convolution}
Bernig, Andreas; Fu, Joseph H. G. Convolution of convex valuations.
Geom. Dedicata 123 (2006), 153--169.
\bibitem{bernig-fu-annals}
Bernig, Andreas; Fu, Joseph H. G.; Hermitian integral geometry. Ann.
of Math., to appear.

\bibitem{chern-45}
Chern, Shiing-shen; On the curvatura integra in a Riemannian
manifold. Ann. of Math. (2) 46, (1945). 674–-684.

\bibitem{federer-curv}
Federer, Herbert; Curvature measures. Trans. Amer. Math. Soc. 93
1959 418–491.

\bibitem{fu-1}
Fu, Joseph H. G.; Curvature measures and generalized Morse theory,
J. Differential Geom. 30 (1989), 619-642.
\bibitem{fu-2}
Fu, Joseph H. G.; Monge-Amp\`ere functions, I, Indiana Univ. Math.
J. 38 (1989), 745-771.
\bibitem{fu-3}
Fu, Joseph H. G.; Monge-Amp\`ere functions, II, Indiana Univ. Math.
J. 38 (1989), 773-789.
\bibitem{fu-90}
Fu, Joseph H. G.; Kinematic formulas in integral geometry. Indiana
Univ. Math. J. 39 (1990), no. 4, 1115–-1154.
\bibitem{fu-94}
Fu, Joseph H. G.; Curvature measures of subanalytic sets. Amer. J.
Math. 116 (1994), no. 4, 819–-880.
\bibitem{fu-unitary}
Fu, Joseph H. G.; Structure of the unitary valuation algebra.  J.
Differential Geom. 72 (2006), no. 3, 509--533.
\bibitem{gelfand-graev-rosu}
Gelfand, I. M.; Graev, M. I.; Ro\c su, R.; The problem of integral
geometry and intertwining operators for a pair of real Grassmannian
manifolds. J. Operator Theory 12 (1984), no. 2, 359--383.
\bibitem{gelfand-graev-schapiro}
Gelfand, I. M.; Graev, M. I.;  \v Sapiro, Z. Ja. Integral geometry
on $k$-dimensional planes. (Russian) Funkcional. Anal. i Prilo\v
zen. 1 1967 15--31.

\bibitem{griffiths-harris}
Griffiths, Phillip; Harris, Joseph; {\itshape Principles of
algebraic geometry.} Reprint of the 1978 original. Wiley Classics
Library. John Wiley \& Sons, Inc., New York, 1994.

\bibitem{guillemin-sternberg}
Guillemin, Victor; Sternberg, Shlomo; {\itshape Geometric
asymptotics.} Mathematical Surveys, No. 14. American Mathematical
Society, Providence, R.I., 1977.

\bibitem{hadwiger-1951}
Hadwiger, Hugo; Translationsinvariante, additive und stetige
Eibereichfunktionale. (German) Publ. Math. Debrecen 2, (1951).
81--94.
\bibitem{hadwiger-book}
Hadwiger, Hugo; {\itshape Vorlesungen \"uber Inhalt, Oberfl\"ache
und Isoperimetrie.} (German) Springer-Verlag,
Berlin-G\"ottingen-Heidelberg 1957.

\bibitem{helgason}
Helgason, Sigurdur; {\itshape The Radon transform.} Second edition.
Progress in Mathematics, 5. Birkh\" auser Boston, Inc., Boston, MA,
1999.

\bibitem{hormander-book}
H\"ormander, Lars; {\itshape The analysis of linear partial
differential operators. I. Distribution theory and Fourier
analysis.} Reprint of the second (1990) edition [Springer, Berlin;
MR1065993 (91m:35001a)]. Classics in Mathematics. Springer-Verlag,
Berlin, 2003.

\bibitem{kashiwara-schapira}
Kashiwara, Masaki; Schapira, Pierre; {\itshape Sheaves on
manifolds.} With a chapter in French by Christian Houzel.
Grundlehren der Mathematischen Wissenschaften [Fundamental
Principles of Mathematical Sciences], 292. Springer-Verlag, Berlin,
1990.

\bibitem{klain-hadwiger}
Klain, Daniel A.; A short proof of Hadwiger's characterization
theorem. Mathematika 42 (1995), no. 2, 329--339.
\bibitem{klain-even}
Klain, Daniel A. Even valuations on convex bodies. Trans. Amer.
Math. Soc. 352 (2000), no. 1, 71--93.
\bibitem{klain-rota}
Klain, Daniel A.; Rota, Gian-Carlo; {\itshape Introduction to
geometric probability.} Lezioni Lincee. [Lincei Lectures] Cambridge
University Press, Cambridge, 1997.

\bibitem{khovanskii-pukhlikov}
Khovanskii, Askold G.; Pukhlikov, Alexander V.; Finitely additive
measures of virtual polyhedra. (Russian) Algebra i Analiz 4 (1992),
no. 2, 161--185; translation in St. Petersburg Math. J. 4 (1993),
no. 2, 337--356.


\bibitem{ludwig-adv-02}
Ludwig, Monika; Projection bodies and valuations. Adv. Math. 172
(2002), no. 2, 158--168.

\bibitem{ludwig-duke-03}
Ludwig, Monika; Ellipsoids and matrix-valued valuations. Duke Math.
J. 119 (2003), no. 1, 159--188.

\bibitem{ludwig-ajm-06}
Ludwig, Monika; Intersection bodies and valuations. Amer. J. Math.
128 (2006), no. 6, 1409--1428.

\bibitem{ludwig-reitzner-annals}
Ludwig, Monika; Reitzner, Matthias; A classification of $SL(n)$
invariant valuations. Ann. of Math., to appear.
\bibitem{ludwig-reitzner-99}
Ludwig, Monika; Reitzner, Matthias; A characterization of affine
surface area. Adv. Math. 147 (1999), no. 1, 138--172.

\bibitem{mcmullen-77}
McMullen, Peter; Valuations and Euler-type relations on certain
classes of convex polytopes. Proc. London Math. Soc. (3) 35 (1977),
no. 1, 113--135.
\bibitem{mcmullen-80}
McMullen, Peter; Continuous translation-invariant valuations on the
space of compact convex sets. Arch. Math. (Basel) 34 (1980), no. 4,
377--384.

\bibitem{mcmullen-survey}
McMullen, Peter; Valuations and dissections. {\itshape Handbook of
convex geometry}, Vol. A, B, 933--988, North-Holland, Amsterdam,
1993.

\bibitem{mcmullen-schneider}
McMullen, Peter; Schneider, Rolf; Valuations on convex bodies.
{\itshape Convexity and its applications}, 170--247, Birkh\"auser,
Basel, 1983.

\bibitem{santalo-book}
Santal\'o, Luis A.; {\itshape Integral geometry and geometric
probability.} With a foreword by Mark Kac. Encyclopedia of
Mathematics and its Applications, Vol. 1. Addison-Wesley Publishing
Co., Reading, Mass.-London-Amsterdam, 1976.

\bibitem{schapira}
Schapira, Pierre; Tomography of constructible functions. (English
summary) Applied algebra, algebraic algorithms and error-correcting
codes (Paris, 1995), 427–435, Lecture Notes in Comput. Sci., 948,
Springer, Berlin, 1995.

\bibitem{schneider-book}
Schneider, Rolf; {\itshape Convex bodies: the Brunn-Minkowski
theory.} Encyclopedia of Mathematics and its Applications, 44.
Cambridge University Press, Cambridge, 1993.
\bibitem{schneider-simple}
Schneider, Rolf; Simple valuations on convex bodies. Mathematika 43
(1996), no. 1, 32--39.
\bibitem{schneider-schuster}
Schneider, Rolf; Schuster, Franz E.; Rotation equivariant Minkowski
valuations. Int. Math. Res. Not. 2006, Art. ID 72894.

\bibitem{viro}
Viro, Oleg; Some integral calculus based on Euler characteristic.
{\itshape Topology and geometry---Rohlin Seminar}, 127--138, Lecture
Notes in Math., 1346, Springer, Berlin, 1988.
\bibitem{wintgen}
P. Wintgen; Normal cycle and integral curvature for polyhedra in
riemannian manifolds, Differential Geometry (G. Soos and J. Szenthe,
eds.), North-Holland, Amsterdam, 1982.
\bibitem{zahle}
Z\"ahle, Martina; Approximation and characterization of generalised
Lipschitz-Killing curvatures. Ann. Global Anal. Geom. 8 (1990), no.
3, 249–-260.












\end{thebibliography}
\end{document}